\UseRawInputEncoding
\documentclass[11pt]{amsart}
\usepackage{bm}
\usepackage{bbm}
\usepackage{nicefrac}\usepackage{amssymb,amstext,amsmath,amsfonts,amsthm}\usepackage{latexsym}
\usepackage{dsfont}
\usepackage{url}
\usepackage{hyperref}
\usepackage{geometry}
\usepackage{color}
\usepackage{anysize}
\marginsize{2.0cm}{2.0cm}{2.0cm}{2.0cm}
\newtheorem{theorem}{Theorem}[section]
\newtheorem{lemma}[theorem]{Lemma}
\newtheorem{corollary}[theorem]{Corollary}
\newtheorem{definition}[theorem]{Definition}
\newtheorem{remark}[theorem]{Remark}
\newcommand{\B}{\textrm{B}}
\newcommand{\N}{\mathbb{N}}
\newcommand{\Z}{\mathbb{Z}}
\newcommand{\R}{\mathbb{R}}
\newcommand{\C}{\mathbb{C}}

\newcommand{\E}[1]{\mathbb{E}\left(#1\right)}
\newcommand{\diag}[1]{\textnormal{diag}\left(#1\right)}
\newcommand{\rL}{\right\}}
\newcommand{\lL}{\left\{}
\newcommand{\rC}{\right ]}
\newcommand{\lC}{\left [}
\newcommand{\rP}{\right)}
\newcommand{\lP}{\left(}
\newcommand{\rabs}{\right|}
\newcommand{\labs}{\left|}
\newcommand{\abs}[1]{\left|�#1 \right|}
\newcommand{\ud}{\mathrm{d}}
\newcommand{\Qr}{\mathbb{P}}
\newcommand{\Prob}[1]{\mathbb{P}\left\{#1\right\}}

\title{Salem--Zygmund inequality for locally sub-Gaussian random variables, random trigonometric polynomials, and random circulant matrices}
\author{Gerardo Barrera}
\address{University of Helsinki, Department of Mathematical and Statistical Sciences.
Exactum in Kumpula Campus, Pietari Kalmin katu 5.
Postal Code: 00560. Helsinki, Finland.}
\email{gerardo.barreravargas@helsinki.fi}
\thanks{*Corresponding author: Gerardo Barrera. Email: gerardo.barreravargas@helsinki.fi}
\thanks{G. Barrera was supported by CIMAT and PIMS}
\author{Paulo Manrique}
\address{National Polytechnic Institute.
Postal code: 07738. Mexico city, Mexico.}
\email{pmanriquem@ipn.mx}
\thanks{P. Manrique was supported by C\'atedras-CONACyT, M\'exico.}
\date{\today}
\begin{document}
\begin{abstract}
In this manuscript we give an extension of the classic Salem--Zygmund inequality for locally sub-Gaussian random variables. As an application, the concentration of the roots of a Kac polynomial is studied, which is the main contribution of this manuscript. 
More precisely, we assume the existence of the moment generating function for the iid random coefficients for the Kac polynomial and prove that there exists an annulus of width
\[\textnormal{O}( n^{-2}(\log n)^{-1/2-\gamma}), \quad \gamma>1/2\] around the unit circle that does not contain roots with high probability. As an another application, we show that the smallest singular value of a random circulant matrix is at least $n^{-\rho}$, $\rho\in(0,1/4)$ with probability $1-\textnormal{O}( n^{-2\rho})$.
\end{abstract}
\maketitle
\section{Introduction}
\markboth{Salem--Zygmund inequality}{Salem--Zygmund inequality}
A classical problem in  Harmonic Analysis is the {\it{quantification}} of the magnitude of the modulus for a  trigonometric polynomial on the unit circle. Erd\"{o}s \cite{erdos} studied the trigonometric polynomial $T_n(x)=\sum_{j=0}^{n-1} \alpha_j e^{ijx}$, $x\in [0,2\pi]$, for choices of signs $\pm 1$ for all $\alpha_j$, and estimated how large $\abs{T_n(x)}$ for $x\in[0,2\pi)$ can be. Salem and Zygmund \cite{sazy} proved that almost all choices of signs satisfy 
\begin{equation}\label{eqn190320181758}
c_1\left(n\log n\right)^{\nicefrac{1}{2}}
\leq
\max\limits_{x\in [0,2\pi]} \left|T_n(x)\right|
\leq 
c_2
\left( n\log n\right)^{\nicefrac{1}{2}}\quad
\textrm{for some positive constants } c_1\;\textrm{ and} \; c_2.
\end{equation}
Inequalities of type \eqref{eqn190320181758} are known as Salem--Zygmund inequality. 
There are different versions of Salem--Zygmund inequality that appear in many areas of modern analysis, see \cite{DM}. In a probabilistic context, the common version of Salem--Zygmund inequality is usually established when the coefficients $\alpha_1,\ldots,\alpha_{n-1}$ of $T_n$ are iid sub-Gaussian random variables, see Chapter 6 in \cite{Kah1985}. 
In the present manuscript, we 
give an extension of Salem--Zygmund inequality for {\it{locally sub-Gaussian}} random coefficients.
This extension allows us to 
study the localization of the roots of a random {\it Kac polynomial} and the probability for the  singularity of a {\it random circulant matrix}.
\subsection{Roots of random trigonometric polynomials} The study of the roots of a polynomial is an old topic in Mathematics. There are formulas to compute the roots for polynomials of degree 2, degree 3 (Tartaglia--Cardano's formula), degree 4 (Ferrari's formula), but due to  Galois' work, for a generic polynomial of degree 5 or more it is not possible to find explicit formulas for computing its roots in terms of radicals. 

For a random polynomial, Bloch and Polya \cite{BP} considered a random polynomial with  iid Rademacher random variables (uniform distribution on $\{-1,1\}$) and proved that the expected number of real zeros are $\mbox{O}({n}^{\nicefrac{1}{2}})$. 
In a series of papers between 1938 and 1939, Littlewood and Offord gave a better bound for the number of real roots of  a random polynomial with iid random coefficients for the cases of Rademacher, Uniform$[-1,1]$, and standard Gaussian \cite{Lub2016}. Kac \cite{kac1943} established his famous integral formula for the density of the number of real roots of a random polynomial with iid coefficients with standard Gaussian distribution. Those were the first steps in the study of roots of random functions, which nowadays is a relevant part of modern Probability and Analysis. For further details, see \cite{DNV} and the references therein.

The localization of the roots of a polynomial
is in general a hard problem. However, there are relevant results in the theory of random polynomials \cite{BhaSam2014}. 
For instance, 
for iid non-degenerate random coefficients with finite logarithm moment, the roots cluster asymptotically near the unit circle and the arguments of the roots are asymptotically uniform distributed. More precisely, Ibragimov and Zaporozhets \cite{IbrZap2013} showed that for a {\it{Kac polynomial}}
\begin{equation}\label{kacp}
G_n(z)=\sum_{j=0}^{n-1} \xi_j z^j \quad
\textrm{ for } z\in \mathbb{C}, 
\end{equation} 
with (real or complex) iid non-degenerate coefficients satisfying $\E{\log(1+|\xi_0|)}<\infty$, its roots are concentrated around the unit circle as $n\to\infty$, almost surely. Moreover, they proved that the condition $\E{\log(1+|\xi_0|)}<\infty$ is necessary and sufficient for the roots of $G_n$ to be asymptotically near the unit circle. 

For iid standard Gaussian random coefficients of $G_n$,
most of the roots are concentrated in an annulus of width $\nicefrac{1}{n}$ centered in the unit circle. However, the nearest root to the unit circle is at least a distance $\mbox{O}({n^{-2}})$, for further details see \cite{Mezin}. Larry and Vanderbei \cite{larry} conjectured that the last statement holds not only for standard Gaussian coefficients but also for Rademacher coefficients. This conjecture was proved by Konyagin and Schlag \cite{Kon1999}. Our 
Theorem~\ref{thm0108201602} establishes that most of the roots of $G_n$ are near to the unit circle in a distance at least $\mbox{O}({n^{-2}} \lP\log n\rP^{-1/2-\gamma})$ for $\gamma>\nicefrac{1}{2}$ with probability $1-\textnormal{O}((\log n)^{-\gamma+1/2})$. 
Konyagin and Schlag \cite{Kon1999} showed that if $G_n$ has iid Rademacher or standard Gaussian random coefficients, then for all $\varepsilon>0$ and large $n$ the following expression \[\min_{z\in\C : \abs{|z|-1}<\varepsilon n^{-2}}|G_n(z)|\geq\varepsilon n^{-\nicefrac{1}{2}}\] holds with probability at least $1-C\varepsilon$, for some positive constant $C$. Karapetyan \cite{Kar1997, Kar1998} studied the sub-Gaussian case, but up to our knowledge, his proof is not complete. Even so, using our extension of Salem--Zygmund inequality and the notion of {\it least common denominator}, which was developed to study the singularity of the random matrices \cite{RV1}, we show that for fixed $t\geq 1$,
\[
\min_{z\in\C\; :\; \abs{\abs{z}-1} \leq  tn^{-2}\lP\log n \rP^{-1/2-\gamma}} \left|G_n(z)\right| \geq t n^{-\nicefrac{1}{2}}(\log n)^{-\gamma},
\] with probability at least $1-\textnormal{O}((\log n)^{-\gamma+1/2})$.
The techniques using in the present paper are not the same using in Konyagin and Schlag \cite{Kon1999}.
The  main result of Konyagin and Schlag only holds for Rademacher and Gaussian iid random coefficients. They did a refined analysis of the characteristic function and applying the so-called circle method. This approach is not straightforward to apply for more general random coefficients, even sub-Gaussian or with finite moment generating function (mgf for short).

The novelty of  this manuscript is the use of the notion of least common denominator for cover more general random coefficients.
This approach works for quite general random coefficients. However, the authors still working for relaxing the assumption of the existence of a mgf.
The main obstacle for relaxing this assumption arises in the control of the maximum modulus over the unit circle of the random polynomial under
the assumption of the existence of some $p-$moment.
We emphasize that the proof is not a direct consequence of \cite{RV2} since good estimates of the least common denominator typically are difficult to obtain. 
We remark that this result and the main result in \cite{Kon1999},  up to our knowledge, are not direct consequences of the so-called concentration inequalities.
\subsection{Random circulant matrices} 
Recall, an $n\times n$ complex circulant matrix, 
denoted by $\textnormal{circ}(c_0,\ldots,c_{n-1})$, has the form
\begin{equation*}
\textnormal{circ}(c_0,\ldots,c_{n-1}):=\left[\begin{array}{ccccc}
c_0 & c_{1} & \cdots & c_{n-2} & c_{n-1} \\
c_{n-1} & c_{0} & \cdots & c_{n-3} & c_{n-2} \\
\vdots & \vdots & \ddots & \vdots & \vdots \\
c_2 & c_3 & \cdots & c_0 & c_1\\
c_{1} & c_{2} & \cdots & c_{n-1} & c_{0}
\end{array}\right],
\end{equation*}
where $c_0,\ldots,c_{n-1}\in \mathbb{C}$. For $\xi_0,\ldots,\xi_{n-1}$ being random variables, we say that 
\[\mathcal{C}_n:=\textnormal{circ}(\xi_0,\ldots \xi_{n-1})\] is an $n\times n$ random circulant matrix.
The circulant matrices are a very common object in different areas of mathematics \cite{HP,KHA,RAU}. In particular, circulant matrices play a crucial role in the study of large-dimensional Toeplitz matrices \cite{AruRajKou2009,SenVir2013}.  
In the theory of the random matrices, the singularity is one aspect that has been intensively studied during recent years \cite{BorCha2012,RV,RV1}. In the case of the random circulant matrices have Rademacher entries, Meckes \cite{Mar2009} proved that the probability of a random circulant matrix is singular tends to zero when its dimension is growing.

As a consequence of our concentration result of the roots for Kac polynomials, 
for a random circulant matrix with iid  zero-mean entries and finite mfg, 
it follows that for all fixed $t\geq 1$ and $\gamma>1/2$, the smallest singular value $s_n\lP\mathcal{C}_n\rP$ of $\mathcal{C}_n$ satisfies \[s_n(\mathcal{C}_n)\geq tn^{-1/2}\lP\log n\rP^{-\gamma}\] with probability $1-\textnormal{O}\left((\log n)^{-\gamma+1/2}\right)$.
However, under weaker assumptions (see below the condition \eqref{fanto}), for $\rho\in(0,1/4)$ we also show \[s_n(\mathcal{C}_n)\geq  n^{-\rho}\] with probability $1-\textnormal{O}\lP n^{-2\rho}\rP$.

The manuscript is organized as follows. In Section~\ref{mainR} we state the main results and their consequences. In 
Section~\ref{sec301120181633} we give the proof of a Salem--Zygmund inequality for random variables with mgf. In 
Section~\ref{sec04082016m02} with the help of Salem--Zygmund inequality and the notion of least common denominator we prove 
Theorem~\ref{thm0108201602} about the location of the roots of a Kac polynomial. Finally, in Section~\ref{sec060320191036} we 
prove Theorem~\ref{thm01042018} about
that the smallest singular value of a random circulant is relatively large with high probability. 
\section{Main results}\label{mainR}
\subsection{Salem--Zygmund inequality}
Recall that a real-valued random variable $\xi$ is said to be 
sub-Gaussian if its mgf is bounded by the mgf of a Gaussian random variable, i.e., there is $b> 0$ such that 
\[
\mathbb{E}(e^{t\xi})\leq e^{\nicefrac{b^2t^2}{2}} \quad \textrm{ for any } t\in \mathbb{R}.
\] 
When this condition is satisfied for a particular value of $b>0$, we say that $\xi$ is $b$-sub-Gaussian or sub-Gaussian with parameter $b$. In particular, it
is straigforward to show that the mean of a sub-Gaussian random variable is necessarily equal to zero. For more details see \cite{chareka} and the references therein.

According to \cite{chareka}, a random variable $\xi$ is called {\it{locally sub-Gaussian}} when its mgf $M_{\xi}$ exists in an open interval around zero. Due to this, it is possible to find constants $\alpha\geq 0$, $\delta\in (0,\infty]$ and $\nu\in\R$ such that 
\[M_{\xi}(t) \leq e^{\nu t +\frac{1}{2}\alpha^2 t^2}\quad \textrm{ for any } t\in (-\delta,\delta).
\]
If the mean of $\xi$ is zero and its variance $\sigma^2$ is finite and positive then we can take $\nu=0$ and $\alpha^2>\sigma^2$ for some $\delta>0$ as the next lemma states.

\begin{lemma}[Locally sub-Gaussian r.v.]\label{lema1559}
Let $\xi$ be a random variable such that its mgf $M_\xi$ exists in an interval around zero. Assume that $\E{\xi}=0$ and $\E{\xi^2}=\sigma^2>0$. Then there is a positive constant $\delta$ such that
\[
M_\xi(t) \leq e^{\nicefrac{\alpha^2 t^2}{2}} \quad \textrm{ for any  } t\in (-\delta,\delta)\; 
\mbox{ and }\; \alpha^2 >\sigma^2.
\]
\end{lemma}
The preceding lemma is not suprising, see for instance 
Remark~2.7.9 in \cite{roman}.
Since its proof is simple, we give it here for completeness of the presentation.

\noindent{\bf Proof.} 
Assume that $M_\xi(t)$ is well-defined for any $t\in (-\delta_1,\delta_1)$, for some 
$\delta_1>0$. Then $M_\xi(t)$ has derivatives of all orders at $t=0$. Define $g(t):=e^{\nicefrac{\alpha^2 t^2}{2}}$, for $t\in\R$. Then $g(0)=1$, $g^\prime(0)=0$ and  $g^{\prime\prime}(0)=\alpha^2$. Let $h(t):= g(t) - M_\xi(t)$, for all $t\in (-\delta_1,\delta_1)$. Since $h$ is continuous and $h^{\prime\prime}(0)=\alpha^2-\sigma^2>0$, then there exists $0<\delta<\delta_1$ such that $h^{\prime\prime}(t)>0$, for every $t\in (-\delta,\delta)$. 
Therefore, the function $h$ is convex in the interval $(-\delta,\delta)$.
As $h^\prime(0)=0$ then $0$ is a local minimum of $h$.
Therefore, it follows that $h(t)\geq h(0)=0$, for every $t\in (-\delta,\delta)$. Thus, the result follows. \hfill$\Box$\\

The classic Salem--Zygmund inequality is usually established for iid sub-Gaussian random variables. But thanks to 
Lemma~\ref{lema1559} we are able to extend it to iid locally sub-Gaussian random variables as it is stated in 
Theorem~\ref{lem06082016m02}. 
Even though, Theorem~\ref{lem06082016m02} is interesting on its own,
we stress that it is also crucial for our approach using in the proof of the main result Theorem~\ref{thm0108201602}.

Before presenting Theorem~\ref{lem06082016m02}, we introduce some useful notations. For simplicity, we keep the same notation between the Euclidean norm and the modulus for the complex numbers. Denote by $\mathbb{T}$ the unit circle $\R/(2\pi\Z)$. For any bounded function $f:\mathbb{T}\to \mathbb{C}$, the infinite norm of $f$ is defined as $\|f\|_\infty=\sup\limits_{x\in\mathbb{T}}|f(x)|$, and $\stackrel{\mathcal{D}}{=}$ means ``equal in distribution''.
\begin{theorem}[Salem--Zygmund inequality for locally sub-Gaussian random variables]\label{lem06082016m02}
Let $\xi$ be a  random variable 
with zero mean and finite positive variance.
Assume that the mgf $M_\xi$ of $\xi$ exists in an open interval around zero. 
Let $\{\xi_k:k\geq 0\}$ be a sequence of iid random variables with $\xi_k\stackrel{\mathcal{D}}{=}\xi$ for every $k\geq 0$.
Let $\phi:[0,1] \rightarrow \mathbb{R}$ is a non-zero continuous function. 
Consider 
$W_n(x)=\sum_{j=0}^{n-1} \xi_j \phi(\nicefrac{j}{n})e^{ijx}$ for any $x\in \mathbb{T}$.
Then, for all large $n$
\[
\mathbb{P}\Big( \|W_n\|_\infty \geq C_0(\lP\log n\rP\sum_{j=0}^{n-1} |\phi(\nicefrac{j}{n})|^2 )^{\nicefrac{1}{2}}\Big)\leq \frac{C_1}{n^2},
\]
where $C_0$ and $C_1$ are positive constants that only depend on the mgf of $\xi$ and the function $\phi$.
\end{theorem}

Actually, under the assumption of finite second moment, a version of a Salem--Zygmund type inequality can be obtained in terms of the expected value of the infinite norm of a random trigonometric polynomial, for more details see \cite{weber}. 
Theorem~\ref{lem06082016m02} provides  an upper bound of how large the infinite norm of a random trigonometric polynomial is in probability. Moreover, Theorem~\ref{lem06082016m02} gives a better bound than Corollary~2 in \cite{weber} as we see below.

Let $\{\xi_k:k\geq 0\}$ be a sequence of iid random variables such that $\E{\xi_0}=0$ and $\E{\xi_0^2}=\sigma^2>0$.  By Corollary~2 in \cite{weber} we have
\begin{align*}
\mathbb{E}\left(\max\limits_{x\in \mathbb{T}}\left|\sum_{j=0}^{n-1} \xi_j e^{ijx}\right|\right)\leq & \; C
\min\left\{(n\log(n+1))\mathbb{E}{(|\xi_0|^2)})^{\nicefrac{1}{2}},n\mathbb{E}{|\xi_0|}\right\} \\
\leq &\; C (n\log(n+1))\mathbb{E}{(|\xi_0|^2)})^{\nicefrac{1}{2}},
\end{align*}
where $C$ is a universal positive constant. By the Markov inequality we obtain
\[
\Qr
\left(\max\limits_{x\in \mathbb{T}}\left|\sum_{j=0}^{n-1} \xi_j e^{ijx}\right| \geq C_0\left(n\log n\right)^{\nicefrac{1}{2}} \right)\leq 
\frac{C (n\log(n+1))\mathbb{E}{(|\xi_0|^2)})^{\nicefrac{1}{2}}}{C_0\left(n\log n\right)^{\nicefrac{1}{2}}}.
\]
Note that the upper bound asymptotically equals a positive constant. On the other hand, under the assumptions of 
Theorem~\ref{lem06082016m02} we deduce
\[
\Qr\left( \max\limits_{x\in \mathbb{T}}\left|\sum_{j=0}^{n-1} \xi_j e^{ijx}\right| \leq C_0\left(n\log n\right)^{\nicefrac{1}{2}} \right)\geq 1-\frac{C_1}{n^2}
\]
for all large $n$, where $C_0$ and $C_1$ are positive constants that only depend on the mgf of $\xi_0$.
\subsection{Kac polynomials}
For using the concept of {\it least common denominator} we introduce the following condition.
We say that a random variable $\xi_0$ satisfies the condition  \eqref{fanto} if
\begin{equation}\label{fanto}
\tag{H}
\sup_{u\in\mathbb{R}} \Prob{ |\xi_0 - u| \leq 1} \leq 1-q\;\; \mbox{ and } \;\;\Prob{ |\xi_0|>M} \leq q/2 
\quad 
\textrm{ for some }  M>0\; \textrm{ and }\;q\in (0,1).
\end{equation}
The notion of concentration function was introduced by P.  L\'evy  in the
context of the study of distributions of sums of random variables.
For $\xi_0$ being  not degenerate, zero mean with mgf, one can deduce that 
condition \eqref{fanto} is valid for some $M>0$ and $q\in (0,1)$. We refer to 
\cite{tikomirov2016}.

The  main result of this manuscript is the following theorem.
\begin{theorem}\label{thm0108201602}
Let $\xi$ be a  random variable 
with zero mean and finite positive variance.
Assume that the mgf $M_\xi$ of $\xi$ exists in an open interval around zero.
Let $\{\xi_k:k\geq 0\}$ be a sequence of iid random variables with $\xi_k\stackrel{\mathcal{D}}{=}\xi$ for every $k\geq 0$. Let
\[ 
\mathcal{M}_n := \lL \min_{z\in\C\; :\; \abs{\abs{z}-1} \leq  tn^{-2}\lP\log n\rP^{-\nicefrac{1}{2}-\gamma}} \left|G_n(z) \right|\leq t n^{-\nicefrac{1}{2}}(\log n)^{-\gamma} \rL.
\]
Then for any fixed $t\geq  1$,
\[
\Qr\left( \mathcal{M}_n \right) = \textnormal{O}\left((\log n)^{-\gamma+1/2}\right),\] where
$\gamma>\nicefrac{1}{2}$ and the implicit constant in the O-notation depends on $t$ and the mgf of $\xi$.
\end{theorem}
\begin{remark}
Observe that all bounded random variables satisfy \eqref{fanto} in Theorem~\ref{thm0108201602} (with a suitable scaling). In particular, the Rademacher distribution which corresponds to the uniform distribution on $\{-1,1\}$ and the uniform distribution on  the interval $[-1,1]$ satisfy \eqref{fanto}. 
\end{remark}
\subsection{Random circulant matrices}
It is well-known that any circulant matrix can be diagonalized in $\mathbb{C}$ using a Fourier basis.
Indeed, 
let $\omega_n:=\exp\left(i\frac{2\pi}{n}\right)$, $i^2=-1$, and $F_n=\frac{1}{\sqrt{n}}(\omega^{jk}_n)_{0\leq j,k\leq n-1}$.
The matrix $F_n$ is called the {\it Fourier matrix} of order $n$. Note that $F_n$ is a unitary matrix. By a straightforward computation it follows
\[\textnormal{circ}(c_0,\ldots,c_{n-1}) = F^*_n \diag{G_n(1),G_n(\omega_n),\ldots,G_n(\omega_n^{n-1})} F_n,\]
where
$G_n$ is the polynomial given by
$G_n(z):=\sum_{k=0}^{n-1}c_kz^k$. Hence, the eigenvalues of $\textnormal{circ}(c_0,\ldots,c_{n-1})$ are $G_n(1),G_n(\omega_n),
\ldots,G_n(\omega_n^{n-1}),$ or equivalently 
\begin{equation}\label{eqn401}
G_n(\omega_n^k)=\sum\limits_{j=0}^{n-1} c_j\exp\left(i\frac{2\pi kj}{n}\right)\quad \textrm{ for any} \quad k=0,\ldots,n-1.
\end{equation}
Expressions like \eqref{eqn401} appear naturally in the study of Fourier transform of periodic functions. For a complete understanding of circulant matrices, we recommend the monograph \cite{David2012}.

In the sequel, we consider an $n\times n$ random circulant matrix $\mathcal{C}_n$, i.e.,
$
\mathcal{C}_n:=\textnormal{circ}(\xi_0,\ldots,\xi_{n-1})$,
where $\xi_0,\ldots,\xi_{n-1}$ are independent random variables.
The smallest singular value of the random circulant matrix $\mathcal{C}_n$ is given by
\begin{equation}\label{holds}
s_n(\mathcal{C}_n) = \min_{0\leq k\leq n-1} |G_n(\omega_n^k)|.
\end{equation}
We remark that in general the smallest singular value is not equal to the smallest eigenvalue
modulus. Since $\mathcal{C}_n$ is a normal matrix, its singular values are the modulus of its eigenvalues. Thus, the following corollary is a direct consequence of Theorem~\ref{thm0108201602}.
\begin{corollary} \label{thm04082016m01}
Let $\xi$ be a  random variable 
with zero mean and finite positive variance.
Assume that the mgf $M_\xi$ of $\xi$ exists in an open interval around zero. 
Let $\{\xi_k:k\geq 0\}$ be a sequence of iid random variables with $\xi_k\stackrel{\mathcal{D}}{=}\xi$ for every $k\geq 0$.
Let $\mathcal{C}_n:=\textnormal{circ}(\xi_0,\ldots,\xi_{n-1})$ be an $n\times n$ random circulant matrix and 
let $s_n(\mathcal{C}_n)$ be the smallest singular value of $\mathcal{C}_n$.
Then, for all fixed $t\geq 1$ and $\gamma > \nicefrac{1}{2}$ we have
\begin{equation}\label{upt}
\Qr\left(s_n(\mathcal{C}_n)\leq tn^{-1/2}\lP\log n\rP^{-\gamma}\right) = \textnormal{O}\lP \lP\log n\rP^{-\gamma+\nicefrac{1}{2}}\rP.
\end{equation}
\end{corollary}
It is possible to weaken the assumptions of 
Corollary~\ref{thm04082016m01}. Using similar reasoning as in the proof of Theorem~\ref{thm0108201602} we obtain the following theorem.
\begin{theorem}\label{thm01042018}
Let $\xi$ be a non-degenerate random variable which satisfies  \eqref{fanto}. Let $\{\xi_k:k\geq 0\}$ be a sequence of iid random variables with $\xi_k\stackrel{\mathcal{D}}{=}\xi$ for every $k\geq 0$.
Let $\mathcal{C}_n:=\textnormal{circ}(\xi_0,\ldots,\xi_{n-1})$ be an $n\times n$ random circulant matrix. Then, for each $\rho\in(0,1/4)$ we have
\begin{equation*}\label{upt2}
\Qr\left(s_n(\mathcal{C}_n)\leq  n^{-\rho}\right) = \textnormal{O}({n^{-2\rho}}).
\end{equation*}
\end{theorem}
\section{Proof of Theorem~\ref{lem06082016m02}. Salem--Zygmund inequality for locally sub-Gaussian random variables}
\label{sec301120181633} 
Firstly, we provide the proof of the following claim which is an important fact that we use in the proof of 
Theorem~\ref{lem06082016m02}.

\noindent
{\bf{Claim ${\bf{1}}$}:}
There exists a random interval $I\subset \mathbb{T}$ of length $\nicefrac{1}{\rho_n}$ with $\rho_n=\nicefrac{3n}{8}$ such that
\begin{align*}\label{claim}
|W_n(x)|\geq \frac{1}{2}
\|W_n\|_\infty\quad \textrm{ for any } x\in I.
\end{align*}
\noindent{\bf Proof.}  Let $p_n(x):=\sum_{j=0}^{n-1} b_j e^{ijx}$, $x\in \mathbb{T}$ be a trigonometric polynomial on $\mathbb{T}$, where
$b_0,\ldots ,b_{n-1}$  are real numbers.
For $x\in \mathbb{T}$ write 
\begin{equation}\label{eq:g}
g_n(x):=|p_n(x)|^2=\left(\sum_{j=0}^{n-1} b_j \cos(jx) \right)^2 + \left(\sum_{j=0}^{n-1} b_j \sin(jx) \right)^2 
\end{equation}
and
\[
h_n(x):= \left( \sum_{j=0}^{n-1} jb_j \cos(jx) \right)^2 + \left(\sum_{j=0}^{n-1} jb_j \sin(jx) \right)^2.
\]
Then
\begin{equation}\label{eq:pp}
\|p_n\|^2_\infty= \sup_{x\in\mathbb{T}} g_n(x)=\|g_n\|_\infty
\quad \textrm{ and } \quad
\|p^{\prime}_n\|^2_\infty=\sup_{x\in\mathbb{T}} h_n(x).
\end{equation}
Recall the Bernstein inequality  $\|p^{\prime}_n\|_\infty\leq n\|p_n\|_\infty$ (see for instance Theorem~14.1.1, Chapter~14, 
page~508 in \cite{RahSch2002}).
For any $x\in \mathbb{T}$ we have
\begin{equation}\label{paulo}
\labs g^{\prime}_n(x) \rabs  \leq  4 \|p_n\|_\infty \|p^{\prime}_n\|_\infty  \leq 4n \|p_n\|^2_\infty = 4n \|g_n\|_\infty.
\end{equation}
Since $g$ is continuous, there exists $x_0\in \mathbb{T}$ such that 
$g(x_0)=\|g_n\|_\infty$. Moreover, by the Mean Value Theorem and relation \eqref{paulo} we obtain
\[
\labs g(x) - g(x_0)\rabs \leq 
\|g^\prime_n\|_\infty \labs x - x_0\rabs \leq 
4n \|g_n\|_\infty \labs x - x_0\rabs
\]
for any $x\in \mathbb{T}$.
Take  $I:=[x_0 - \frac{3}{16n},x_0+\frac{3}{16n}]\subset 
\mathbb{T}$. Notice that the length of $I$ is $\frac{3}{8 n}$. The preceding inequality yields
\[
\labs g(x) - g(x_0)\rabs \leq \frac{3}{4}  \|g_n\|_\infty
\quad \textrm{ for any } x\in I. 
\] 
Since $g(x_0)=\|g_n\|_\infty$, the triangle inequality yields
$(\nicefrac{1}{4}) \|g_n\|_\infty \leq  \labs g_n(x)\rabs$ for any $x\in I$. The preceding inequality with the help of relation \eqref{eq:g} and relation \eqref{eq:pp} implies
\[
\frac{1}{2} \|p_n\|_\infty \leq  \labs p_n(x)\rabs
\;\; \textrm{ for any } x\in I.
\] \hfill$\Box$

Now, we are ready to provide the proof of 
Theorem~\ref{lem06082016m02}.
\begin{proof}[Proof of Theorem~\ref{lem06082016m02}]
By Lemma~\ref{lema1559}, there exists a $\delta >0$ such that
\[
M_\xi(t) \leq e^{\alpha^2 t^2/2}\quad \mbox{ for any  } t\in (-\delta,\delta),\;
\mbox{ where  } \alpha^2 >\sigma^2>0.
\]
For each $j\in \{0,\ldots,n-1\}$, define $f_j(x)=\phi(\nicefrac{j}{n})e^{ijx}$, $x\in\mathbb{T}$.
Let $r_n:=\sum_{j=0}^{n-1} |\phi(\nicefrac{j}{n})|^2$.
At first, we suppose that the $f_j$ are real (we consider only the real part or the imaginary part) and we write $S_n:=\|W_n\|_\infty$. Since $\|f_j\|_\infty\leq \|\phi\|_\infty=:K$ for every $j=0,\ldots,n-1$, we obtain
\begin{align*}
e^{\nicefrac{\alpha^2 t^2 r_n}{2}}&=
\prod_{j=0}^{n-1} e^{\nicefrac{\alpha^2 t^2\|f_j\|^2_\infty}{2}}
\geq \prod_{j=0}^{n-1} e^{\nicefrac{\alpha^2 t^2|f_j(x)|^2}{2}}\geq \prod_{j=0}^{n-1} \E{ e^{t \xi_j f_j(x)}}\\
&=\E{\prod_{j=0}^{n-1} e^{t \xi_j f_j(x)}}=\E{e^{t W_n(x)}}\quad \textrm{ for any } t\in (-\nicefrac{\delta}{K},\nicefrac{\delta}{K}).
\end{align*}
By {\bf{Claim ${\bf{1}}$}},
there exists a random interval $I\subset \mathbb{T}$ of length $\nicefrac{1}{\rho_n}$ with $\rho_n=\nicefrac{8n}{3}$
such that
$W_n(x)\geq  \nicefrac{S_n}{2}$ or $-W_n(x) \geq \nicefrac{S_n}{2}$ on $I$.
Denote
by $\mu$ the normalized Lebesgue measure on $\mathbb{T}$.
Observe that
\begin{align*}
e^{\nicefrac{tS_n}{2}}=\frac{1}{\mu(I)}\int\limits_{I}e^{\nicefrac{tS_n}{2}}\ud x
\leq  \frac{1}{\mu(I)}\int\limits_{I}
\left( e^{tW_n(x)} + e^{-tW_n(x)} \right)
\ud x.
\end{align*}
Then, for every $t\in (-\nicefrac{\delta}{K},\nicefrac{\delta}{K})$ we have
\begin{align*}
\E{e^{\nicefrac{t S_n}{2}}} & \leq  \rho_n \E{ \int_{I} \left( e^{tW_n(x)} + e^{-tW_n(x)} \right) \mu(\ud x) } \\
 & \leq  \rho_n \E{ \int_{\mathbb{T}} \left( e^{tW_n(x)} + e^{-tW_n(x)} \right) \mu(\ud x) } \; \leq \;  2\rho_n e^{\nicefrac{\alpha^2 t^2 r_n}{2}}.
\end{align*} 
The preceding inequality yields
\[
\E{\exp\left\{\frac{t}{2}\left( S_n -\alpha^2 t r_n - \frac{2}{t} \log\left( 2\rho_n l\right)\right)\right\}} \leq \frac{1}{l}
\quad
\textrm {for any } l>0\; \textrm{ and }\;  t\in (-\nicefrac{\delta}{K},\nicefrac{\delta}{K}),
\]
which implies
\[
\Qr\left( S_n \geq \alpha^2 tr_n + \frac{2}{t}\log\left(2\rho_n l \right)\right) \leq \frac{1}{l}\quad
\textrm {for any } l>0\; \textrm{ and }\;  t\in (-\nicefrac{\delta}{K},\nicefrac{\delta}{K}).
\]
Note that $\lim\limits_{n\to \infty}\frac{r_n}{n}=\int_{0}^{1}|\phi(x)|^2\ud x>0$. 
By taking $l_n=cn^2$ where $c$ is a positive constant, 
we have
$\left|\frac{\log(2\rho_n l_n)}{\alpha^2 r_n}\right|<\nicefrac{\delta^2}{K^2}$
for all large $n$. By choosing $t_n=\left(\frac{\log(2\rho_n l_n)}{\alpha^2 r_n}\right)^{\nicefrac{1}{2}}$ we obtain
\[
\Qr\left(S_n \geq 3 \left( \alpha^2 r_n\log\left( 2\rho_n l_n\right)\right)^{\nicefrac{1}{2}}\right) \leq \frac{1}{l_n} \quad \textrm{ for all large } n.
\]
Since $f_j=\mbox{Re}(f_j)+i\mbox{Im}(f_j)$,   we get for all large $n$
\[
\Qr\left( \| \mbox{Re}(W_n) \|_\infty \geq 3 \left( \alpha^2\sum_{j=0}^{n-1} \|\mbox{Re} (f_j)\|^2_\infty\log\left( 2\rho_n l_n\right)\right)^{\nicefrac{1}{2}}\right)  \leq  \frac{1}{l_n}
\]
and
\[
\Qr\left( \| \mbox{Im}(W_n) \|_\infty \geq 3 \left( \alpha^2\sum_{j=0}^{n-1} \|\mbox{Im} (f_j)\|^2_\infty\log\left( 2\rho_n l_n\right)\right)^{\nicefrac{1}{2}}\right)  \leq  \frac{1}{l_n}.
\]
Finally, since $\rho_n=\frac{8n}{3}$, the choose of $l_n=\frac{3n^2}{16}$ yields
\[
\Qr\left( \|W_n\|_\infty \geq 6\alpha \sqrt{3} \left( r_n \log n\right)^{\nicefrac{1}{2}} \right)\leq \frac{32}{3n^2} \quad \textrm{ for all large } n.
\]
\end{proof}
\section{Proof of Theorem~\ref{thm0108201602}. Localization of the roots for Kac polynomials} \label{sec04082016m02}
The proof is based on the small ball probability 
of linear combinations of iid random variables introduced by  Rudelson and Vershynin in \cite{RV2}.
Throughout the proof, $\|\cdot\|_2$ denotes the Euclidean norm, $|\cdot|$ denotes the complex norm and $\det(\cdot)$ the determinant function that acts on the squared matrices. 
We consider the module $\pi$ of a real number $y$,  $y \hspace{-0.2cm}\mod \pi$, which is defined as the set of numbers $x$ such that $x-y =k\pi$ for some $k\in\Z$. 
\begin{definition}[Least common denominator (lcd for short)]
Let $L$ be any positive number and let $V$ be any {\em{deterministic}} matrix of dimension $2\times n$. 
The {\it least common denominator} (lcd) of $V$ is defined as 
\[
D(V):= \inf \lL \|\theta \|_2 >0 : \theta\in\R^2, \;\mathrm{dist}\lP V^T\theta,\Z^n \rP < L\sqrt{\log_{+} \left(\frac{\| V^T\theta\|_2}{L}\right)} \rL,
\]
where 
$\mathrm{dist}(v,\mathbb{Z}^n)$ denotes the distance between the vector $v\in \mathbb{R}^n$ and the set $\mathbb{Z}^n$,
and $\log_{+}:=\max\{\log,0\}$.
\end{definition} For more details of the concept of lcd see Section~7 of \cite{RV2}.
Observe that $D\lP a V\rP = \lP\nicefrac{1}{|a|}\rP D(V)$ for any $a\not=0$. Indeed, from the definition of $D\lP a V\rP$ we have that $D\lP a V\rP\leq \|\theta\|_2$ for any $\theta\in \mathbb{R}^2$ such that
\[
\mbox{dist}\lP (aV)^T\theta,\Z^n \rP < L\sqrt{\log_{+} \left(\frac{\| (aV)^T\theta\|_2}{L}\right)}=L\sqrt{\log_{+} \left(\frac{\| V^T(a\theta)\|_2}{L}\right)}.
\]
Therefore, from the definition of $D(V)$ we deduce $D(V)\leq \|a\theta\|_2=|a|\|\theta\|_2$.
Since $a\not=0$, then $(\nicefrac{1}{|a|})D(V)\leq \|\theta\|_2$.
Again, from the definition of $D(aV)$ we deduce that
$(\nicefrac{1}{|a|})D(V)\leq D(aV)$.
On the other hand, from the definition of $D(V)$ we have that
$D\lP V\rP\leq \|\theta\|_2$ for any $\theta\in \mathbb{R}^2$ such that
\[
\mbox{dist}\lP V^T\theta,\Z^n \rP < L\sqrt{\log_{+} \left(\frac{\| V^T\theta\|_2}{L}\right)}=L\sqrt{\log_{+} \left(\frac{\| (aV)^T(\nicefrac{\theta}{a})\|_2}{L}\right)}.
\]
Therefore, from the definition of $D(aV)$ we deduce $D(aV)\leq \|\nicefrac{\theta}{a}\|_2=\nicefrac{\|\theta\|_2}{|a|}$.
Consequently, $|a|D(aV)\leq \|\theta\|_2$. Again, from the definition of $D(V)$ we deduce that $|a|D(aV)\leq D(V)$. Putting all these pieces together we obtain the next useful lemma.
\begin{lemma}\label{prop05112019} For all $a\neq 0$, the lcd of any matrix $V\in \mathbb{R}^{2\times n}$ satisfies 
$D(V)=|a|D(aV)$.
\end{lemma}
Let $X$ be a {\em{random}} vector  of dimension $n\times 1$ whose entries are iid satisfying \eqref{fanto}. Assume $\det(VV^T)>0$. 
For any $a>0$ and $t\geq 1$, by Theorem~7.5 (Section~7 in \cite{RV2}) we have
\[
\begin{split}\label{rvi}
\Prob{ \|V X\|_2 \leq \frac{t\sqrt{2}}{a} } & =
\Prob{ \|aV X\|_2 \leq \sqrt{2}t } \leq 
\frac{C^2 L^2}{2a^2(\det(VV^T))^{\nicefrac{1}{2}}}
\lP t + \frac{\sqrt{2}}{D(a V)}\rP^2,
\end{split} 
\]
where $L\geq \sqrt{\nicefrac{2}{q}}$ with $q$ given in \eqref{fanto}, $D(aV)$ is the least common denominator of $aV$, and the 
constant $C$ only depends on $M$, $q$.
Recall the well-known inequality  $(x+y)^2\leq 2 x^2+ 2y^2$ for any $x,y\in \mathbb{R}$. By Lemma~\ref{prop05112019}, it follows that $D(aV)= (\nicefrac{1}{a})D(V)$ for all $a>0$. Therefore, 
\begin{align}\label{eqn270220191047}
\Prob{ \|aV X\|_2 \leq \sqrt{2}t } & \leq 
\frac{C^2 L^2}{a^2(\det(VV^T))^{\nicefrac{1}{2}}}t^2+
\frac{2C^2 L^2}{a^2(\det(VV^T))^{\nicefrac{1}{2}}(D(a V))^2} \nonumber\\
& \leq  \frac{C^2 L^2}{a^2(\det(VV^T))^{\nicefrac{1}{2}}}t^2+
\frac{2C^2 L^2}{(\det(VV^T))^{\nicefrac{1}{2}}(D(V))^2}.
\end{align}
In order to obtain a meaningful upper bound for the  left-hand side of the preceding inequality, 
 it is needed to do a refined analysis of the following quantities:
 {\em{a lower bound for $\det(VV^T)$}} and {\em{a lower bound for $D(V)$.}} Implicitly, in the definition of the $D(V)$ we also need to estimate $\| V^T\theta\|_2$ for some adequate $\theta\in \mathbb{R}^2$.

\subsection{\bf{Small ball probability analysis}} 
The following  analysis explains the reason of introducing the concept of the least common denominator, which is a crucial part along the proof of Theorem~\ref{thm0108201602}. Recall
\[
G_n(z)= \sum_{j=0}^{n-1} \xi_j z^j\quad 
\textrm{ for } z\in \mathbb{C}.
\]
For $G_n$, we associate a random trigonometric polynomial
\[
W_n(x)= \sum_{j=0}^{n-1} \xi_j e^{ijx}\quad 
\textrm{ for } x\in \mathbb{T},
\]
where $\mathbb{T}$ denotes the unit circle $\R/(2\pi\Z)$.
Assume $n\geq 2$ and $\gamma>\nicefrac{1}{2}$. Let 
$N=\lfloor n^2 \lP \log n\rP^{\nicefrac{1}{2}+\gamma} \rfloor$ and $x_\alpha = \nicefrac{\alpha}{N}$ for $\alpha\in \{0,1,2,\ldots, N-1\}$.
Let $t\geq 1$  be fixed and let
 $C_0>0$ be the suitable positive constant being given in 
 Theorem~\ref{lem06082016m02}. Define the following event
\[
\mathcal{G}_n : = \lL \| W^{\prime}_n\|_\infty \leq C_0 n^{\nicefrac{3}{2}}\lP \log n\rP^{\nicefrac{1}{2}}, \max_{z\in\C\; :\; \abs{\abs{z}-1} \leq  2tn^{-2}} \abs{G_n(z)} \leq n^{\nicefrac{3}{2}} \rL,
\] 
where $W^{\prime}_n$ denotes the derivative of $W_n$ on $\mathbb{T}$. For short, we also denote by $\Qr\lP A,B\rP$ the probability $\Qr\lP A\cap B\rP$ for any two events $A$ and $B$. Recall
\[ 
\mathcal{M}_n = \lL \min_{z\in\C\; :\; \abs{\abs{z}-1} \leq  tn^{-2}\lP\log n\rP^{-\nicefrac{1}{2}-\gamma}} \left|G_n(z) \right|\leq t n^{-\nicefrac{1}{2}}(\log n)^{-\gamma} \rL.
\]
By the Boole--Bonferroni inequality we obtain
\begin{align}\label{eq:bound}
 \Qr\lP \mathcal{M}_n \rP &  \leq  \Qr\lP \mathcal{M}_n, \mathcal{G}_n\rP + \Qr\lP \| W'_n\|_\infty \geq C_0 n^{\nicefrac{3}{2}}\lP \log n\rP^{\nicefrac{1}{2}} \rP \nonumber\\
&  \quad +\; \Qr\lP \displaystyle\max_{z\in\C\; :\; \abs{\abs{z}-1} \leq  2tn^{-2}}\abs{G_n(z)} \geq n^{\nicefrac{3}{2}} \rP \nonumber\\
& =: \Qr\lP \mathcal{M}_n, \mathcal{G}_n\rP + I_1 + I_2.
\end{align}
Our goal is to show that every probability on the right side of the above expression is decreasing to zero when $n$ tends to infinity.

Using the Berstein inequality (Theorem~14.1.1 in \cite{RahSch2002}) and Theorem~\ref{lem06082016m02} for $\phi \equiv 1$, for all large $n$ we have
\[
\Qr\lP \| W^{\prime}_n\|_\infty \geq C_0 n^{\nicefrac{3}{2}}\lP \log n\rP^{\nicefrac{1}{2}} \rP \leq  \Qr\lP \| W_n\|_\infty \geq C_0 \lP n\log n\rP^{\nicefrac{1}{2}} \rP \leq \frac{C_1}{n^2}.
\]
On the other hand, using the Markov inequality we obtain
\[
\begin{split}
&\Qr\lP \displaystyle\max_{z\in\C\; :\; \abs{\abs{z}-1} \leq  2tn^{-2}}\abs{G_n(z)} \geq n^{3/2} \rP  \leq  \Qr\lP \displaystyle\sum_{j=0}^{n-1} \abs{\xi_j} \lP 1+\frac{2t}{n^2}\rP^{j} \geq n^{3/2}\rP \\
&\hspace{1.5cm} \leq \frac{1}{n^{3/2}} \E{\sum_{j=0}^{n-1} \abs{\xi_j} \lP 1 + \frac{2t}{n^2}\rP^j} \leq \frac{e^{2t}  \E{\abs{\xi_0}}n}{n^{\nicefrac{3}{2}}} = \frac{e^{2t}\E{\abs{\xi_0}}}{n^{\nicefrac{1}{2}}},
\end{split}
\]
where the last inequality follows from the following fact:
for any $j\in \{0,\ldots,n^2\}$ we have
\[
\lP 1 + \frac{2t}{n^2}\rP^{j}\leq \lP 1 + \frac{2t}{n^2}\rP^{n^2}\leq e^{2t}.
\]
Therefore,  
\begin{equation}\label{eqn051120191037}
I_1 + I_2 = \textnormal{O}\lP n^{-1/2}\rP,
\end{equation}  
where the implicit constant depends on the distribution of $\xi_0$ and $t$.
We stress that the rate of the convergence in \eqref{eqn051120191037} can be improved,
 however, the contributed term in the right-hand side of \eqref{eq:bound} is
$\Qr\lP \mathcal{M}_n, \mathcal{G}_n\rP$.

In the sequel, we analyze the strategy to prove that $\mathbb{P}(\mathcal{M}_n, \mathcal{G}_n)$ is small. First, we construct a set of closed balls that covers $\{z\in\C:\abs{\abs{z}-1} \leq  tn^{-2}\}$. For each closed ball, we reduce the event $\{\mathcal{M}_n, \mathcal{G}_n\}$ to a ``simple event'' using Taylor's Theorem. Finally, we use the concept of lcd to show that the probability of each ``simple event'' is sufficiently small.

The strategy is to consider a set of balls centered at a point on the unit circle with a suitable radius. We distinguish two kind of balls. The special balls centered in $1+0i$ and $-1+0i$, where the radius $r$ is {\it large}, $r=2tn^{-11/10}$, and the balls centered in points $z$ with argument satisfying $n^{-11/10}<\abs{\arg(z)\mod \pi } < \pi - n^{-11/10} $ with {\it small} radius, $r=2t n^{-2}$.

Recall that for any $x\in \mathbb{R}$, $\lfloor x \rfloor$  denotes the greatest integer less than or equal to $x$.
Let $N:= \lfloor n^2\lP \log n \rP^{\nicefrac{1}{2}+\gamma} \rfloor$ and $x_\alpha := \frac{\alpha}{N}$ for $\alpha=0,1,\ldots,N-1$. For $a\in\C$ and $s>0$, denote by $\textrm{B}\lP a, s\rP$ the closed ball with center $a$ and radius $s$, i.e., $\textrm{B}\lP a, s\rP = \lL z\in\C : \abs{z-a} \leq  s\rL$. Denote by 
$\mathbb{S}^1$ the unit circle.
Let 
\begin{align*}
\mathcal{A}\left(\mathbb{S}^1,tn^{-2}\lP\log n\rP^{-1/2-\gamma}\right):=\left\{z\in \mathbb{C}:\abs{\abs{z}-1} \leq tn^{-2}\lP\log n\rP^{-1/2-\gamma}\right\}.
\end{align*}
Note that
\begin{align*} 
 \mathcal{A}\left(\mathbb{S}^1,tn^{-2}\lP\log n\rP^{-1/2-\gamma}\right) = & 
\lL z\in\mathcal{A} : n^{-11/10}<\abs{\arg(z) } < \pi - n^{-11/10} \rL \\ 
&\hspace{-3cm}\cup \lL z\in\mathcal{A} : \abs{\arg(z)}\leq n^{-11/10}\quad \mbox{ or }\quad \abs{\arg(z)-\pi} \leq n^{-11/10} \rL.
\end{align*}
Let $t\geq 1$ and observe that
\begin{align*} 
&\lL z\in\mathcal{A} : \abs{\arg(z)}\leq n^{-11/10}\quad \mbox{ or }\quad \abs{\arg(z)-\pi} \leq n^{-11/10} \rL \\
&\hspace{3.5cm} \subset \B\lP -1+0i, 2tn^{-11/10}\rP \cup \B\lP 1+0i, 2tn^{-11/10}\rP.
\end{align*} 
The preceding inclusion yields that any $z\in\mathcal{A}$ with \textit{small argument} belongs in the balls centered at $1+0i$ and $-1+0i$ with radius $2tn^{-11/10}$. On the other hand, for $z\in\mathcal{A}$ with \textit{large argument} we have
\begin{align*}
&\lL z\in\mathcal{A} : n^{-11/10}<\abs{\arg(z) } < \pi - n^{-11/10} \rL \\
&\hspace{1cm}\subset \bigcup^{N-1} _{\substack{\alpha=1 \\ \alpha\;:\; n^{-11/10}<\abs{2\pi x_\alpha\hspace{-0.2cm}\mod\pi} < \pi - n^{-11/10}}}
\B \lP e^{i2\pi x_\alpha},2tn^{-2}\lP\log n\rP^{-\nicefrac{1}{2}-\gamma}\rP.
\end{align*}
We define $[N-1]:=[1,N-1]\cap \mathbb{N}$ and
\[
\begin{split}
J_1(n,N)&:=\left\{\alpha\in [N-1]: \gcd\lP \alpha, N\rP \geq n^{1+1/10} \lP\log n\rP^{-\gamma}\right\},\\
J_2(n,N)&:=\left\{\alpha\in [N-1]: n^{1+1/10} \lP\log n\rP^{-\gamma}\geq \gcd\lP\alpha,N\rP\geq n\lP\log n\rP^{\nicefrac{1}{2}+\gamma}\right\},\\
J_3(n,N)&:=\left\{\alpha\in [N-1]: n\lP\log n\rP^{\nicefrac{1}{2}+\gamma} \geq \gcd\lP \alpha, N\rP \geq n^{9/10}\lP \log n\rP^{\nicefrac{1}{2}+\gamma}\right\},
\end{split}
\]
where $\gcd(\alpha,N)$ denote the greatest common divisor between $\alpha$ and $N$.
Observe that 
for any $\alpha\in J_3(n,N)$ we have
\[
n - \frac{1}{n\lP \log n\rP^{\nicefrac{1}{2}+\gamma}}\leq  \frac{N}{\gcd\lP \alpha, N\rP} \leq n^{11/10}.
\] 
The preceding inequalities yield that the irreducible fraction of $x_\alpha$ is as small as a multiple of $n^{-11/10}$. Therefore,
\begin{align*}
&\bigcup^{N-1} _{\substack{\alpha=1 \\ \alpha\;:\; n^{-11/10}<\abs{2\pi x_\alpha \hspace{-0.2cm}\mod\pi} < \pi - n^{-11/10}}}
\B \lP e^{i2\pi x_\alpha},2tn^{-2}\lP\log n\rP^{-\nicefrac{1}{2}-\gamma}\rP\\
&\hspace{2cm}= \bigcup_{\alpha\in J_1(n,N)}
\B \lP e^{i2\pi x_\alpha},2tn^{-2}\lP\log n\rP^{-\nicefrac{1}{2}-\gamma}\rP\\
&\hspace{2.5cm}
\cup 
\bigcup_{\alpha\in J_2(n,N)}
\B \lP e^{i2\pi x_\alpha},2tn^{-2}\lP\log n\rP^{-\nicefrac{1}{2}-\gamma}\rP\\
&\hspace{2.5cm} \cup 
\bigcup_{\alpha\in J_3(n,N)}
\B \lP e^{i2\pi x_\alpha },2tn^{-2}\lP\log n\rP^{-\nicefrac{1}{2}-\gamma}\rP.
\end{align*}
We emphasize that if $\alpha\in J_1(n,N)\cup J_2(n,N)\cup J_3(n,N)$, then we have
\[n^{-11/10} < \abs{2\pi x_\alpha \hspace{-0.3cm}\mod\pi} < \pi - n^{-11/10}.\] 
Consequently,
\begin{equation}\label{mainestimate}
\begin{split}
\Prob{\mathcal{M}_n, \mathcal{G}_n} \; \leq & \; \Prob{\mathcal{G}_n, \min_{z\in \B\lP 1+0i, 2tn^{-11/10}\rP} \abs{G_n(z)} < tn^{-1/2}\lP \log n\rP^{-\gamma}} \\
&\quad +\Prob{\mathcal{G}_n, \min_{z\in \B\lP -1 + 0i, 2tn^{-11/10}\rP} \abs{G_n(z)} < tn^{-1/2}\lP \log n\rP^{-\gamma}}  \\
&\quad  + \sum\limits_{\alpha\in J_1(n,N)} \Prob{\mathcal{G}_n, \B_\alpha}+\quad
\sum\limits_{\alpha\in J_2(n,N)} \Prob{\mathcal{G}_n, \B_\alpha}\quad+\sum\limits_{\alpha\in J_3(n,N)} \Prob{\mathcal{G}_n, \B_\alpha},
\end{split}
\end{equation}
where
\[
\B_\alpha := \lL \min_{z\in B\lP e^{i2\pi x_\alpha}, 2t n^{-2}\lP \log n\rP^{-\nicefrac{1}{2}-\gamma}\rP} \abs{G_n(z)} < tn^{-1/2}\lP \log n \rP^{-\gamma} \rL.
\]
\subsubsection{\bf Small ball analysis at the points $\bm{1+0i}$ and $\bm{-1+0i}$}\label{subsubsec07052020A} On the two points $1\pm 0i$ we have the largest two closed balls, which are considered in our set of balls. This is remarkable since the number of real roots of a Kac polynomial for some common random variables is at least $\textnormal{O}( \frac{\log n}{\log \log \log n})$ with high probability \cite{Nguyen2016}. This means that the real roots of a Kac Polynomial are moving slowly to the unit circle. 

On the one hand, 
let  $z\in\textnormal{B}\lP 1+0i, 2tn^{-\nicefrac{11}{10}}\rP$. By Taylor's Theorem we obtain
\[
\abs{G_n(z)-G_n(1)}\leq\abs{z-1}\abs{G^{\prime}_n(1)} + \abs{R_2(z)},
\] where $R_2(z)$ is the error of the Taylor approximation of order 2. On the event $\mathcal{G}_n$ we have
\begin{align*}
\abs{R_2(z)} & \leq  \frac{\lP 2tn^{-1-1/10}\rP^{2}}{1-\textnormal{o}(1)} \lC \max_{z\in\C\; :\; \abs{\abs{z}-1} \leq 2tn^{-2}}\abs{G_n(z)}\rC \\
& \leq  \frac{4t^2 n^{-2-1/5} n^{3/2}}{1-\textnormal{o}(1)} = \frac{4t^2 n^{-1/2-1/5}}{1-\textnormal{o}(1)},
\end{align*} where $\textnormal{o}(1) = 2tn^{-1-1/10}$. Assuming that $\mathcal{G}_n$ holds, the preceding inequality yields
\begin{align*}
\abs{G_n(z)-G_n(1)} & \leq  2tn^{-1-1/10} \abs{G_n'(1)} + \frac{4t^2 n^{-1/2-1/5}}{1-\textnormal{o}(1)} \leq  2tn^{-1-1/10} \|W_n'\|_\infty + \frac{4t^2 n^{-1/2-1/5}}{1-\textnormal{o}(1)} \\
& \leq  2C_0 t n^{1/2-1/10} \lP \log n\rP^{1/2} + \frac{4t^2 n^{-1/2-1/5}}{1-\textnormal{o}(1)}.
\end{align*}
Hence,
\begin{eqnarray*}
& \Qr\lP \mathcal{G}_n, \displaystyle\min_{z\in\textnormal{B}\lP 1+0i, 2tn^{-11/10}\rP} \abs{G_n(z)} \leq tn^{-1/2}\lP \log n\rP^{-\gamma}\rP \\
&\hspace{4cm} \leq \Qr\lP\abs{G_n(1)}\leq 2C_2 t n^{1/2-1/10} \lP\log n\rP^{1/2}  \rP,
\end{eqnarray*} where $2C_2 = 2C_0t +4t^2 + 1$. Since $G_n(1)=\sum_{j=0}^{n-1} \xi_j$, Corollary~7.6 in \cite{RV2} implies 
for $L\geq \sqrt{1/q}$ (with $q$ given in \eqref{fanto}) that
\[
\Qr\lL \abs{G_n(1)}\leq 2C_2tn^{1/2-1/10}\lP\log n\rP^{1/2}\rL \leq \frac{C_3 L}{\|\mathbf{a}\|} \lP 2C_2t + \frac{1}{D(\mathbf{a})}\rP,
\] where $C_3$ is a positive constant and $D(\mathbf{a})$ is the lcd of the vector
\[
\mathbf{a}=( n^{1/2-1/10} \lP\log n\rP^{1/2})^{-1} \lP1,\ldots,1\rP \in \R^{n}.
\]
By Proposition 7.4 in \cite{RV2} we have $D(\mathbf{a})\geq \frac{1}{2}n^{1/2-1/10}\lP\log n\rP^{1/2}$. Therefore,
\begin{eqnarray}\label{eqn071120190936}
&\hspace{-3cm}\Qr\lP\abs{G_n(1)}\leq 2C_2 t n^{1/2-1/10} \lP\log n\rP^{1/2}  \rP &
\nonumber\\
&\hspace{2.2cm} \leq  \frac{C_3 L \lP\log n\rP^{1/2}}{n^{1/10}}\lP 2C_2t + \frac{2}{n^{1/2-1/10}\lP \log n\rP^{1/2}}\rP  
\leq  \frac{\lP 2C_2t+2\rP L\lP\log n\rP^{1/2}}{n^{1/10}}. 
\end{eqnarray}

On the other hand, let $z\in\textnormal{B}\lP -1+0i, 2tn^{-11/10}\rP$. Assuming that $G_n$ holds, Taylor's Theorem implies
\begin{align*}
\abs{G_n(z)-G_n(-1)} & \leq  \abs{z+1}\abs{G'_n(-1)} + \abs{R_2(z)} \leq  2tn^{-1-1/10} \|W'_n\|_\infty+ \frac{4t^2n^{-1/2-1/5}}{1-\textnormal{o}(1)} \\
& \leq  \lP2C_0 t + 4t^2\rP n^{1/2-1/10}\lP \log n\rP^{1/2}.
\end{align*} Thus,
\begin{align*}
& \Qr\lP\mathcal{G}_n,\displaystyle\min_{z\in\textnormal{B}\lP -1+0i, 2tn^{-11/10}\rP} \abs{G_n(z)} \leq tn^{-1/2}\lP \log n\rP^{-\gamma}\rP \\
&\hspace{3cm} \leq \Qr\lP\abs{G_n(-1)}\leq 2C_2 t n^{1/2-1/10} \lP\log n\rP^{1/2}\rP.
\end{align*} 
Since $G_n(-1) = \sum_{j=0}^{n-1} \lP -1\rP^{j}\xi_j$, by Corollary~7.6 in \cite{RV2} for $L\geq \sqrt{1/q}$ (with $q$ given in \eqref{fanto}) we obtain
\[
\Qr\lP\abs{G_n(-1)}\leq 2C_2 t n^{1/2-1/10} \lP\log n\rP^{1/2}\rP \leq \frac{C_3L}{\|\mathbf{b}\|}\lP 2C_2t +\frac{1}{D(\mathbf{b})}\rP,
\] where $C_3$ is a positive constant and $D(\mathbf{b})$ is the lcd of the vector
\[
\mathbf{b}=( n^{1/2-1/10}\lP \log n\rP^{1/2})^{-1}\lP1,-1,\ldots,(-1)^{n-1}\rP\in\R^n.
\] 
By Proposition 7.4 in \cite{RV2}, we have $D(\mathbf{b})\geq \frac{1}{2} n^{1/2-1/10} \lP \log n\rP^{1/2}$. Therefore,
\begin{eqnarray}\label{eqn091120190937}
&\hspace{-4cm}\Qr\lP\abs{G_n(-1)}\leq 2C_2 t n^{1/2-1/10} \lP\log n\rP^{1/2}\rP  \nonumber\\
 &\hspace{3.5cm}\leq  \frac{C_3 L \lP \log n\rP^{1/2}}{n^{1/10}} \lP 2C_2 t + \frac{2}{n^{1/2-1/10}\lP \log n\rP^{1/2}}\rP 
\leq  \frac{\lP 2C_2 t + 2\rP L \lP \log n \rP^{1/2}}{n^{1/10}}. 
\end{eqnarray}
Combining \eqref{eqn071120190936} and \eqref{eqn091120190937} we obtain
\begin{align*}\label{eqn091120190940}
&\Qr\lP \mathcal{G}_n, \min_{z\in\textnormal{B}\lP 1+0i, 2tn^{-\nicefrac{11}{10}}\rP} \abs{G_n(z)} \leq tn^{-\nicefrac{1}{2}}\lP \log n\rP^{-\gamma}\rP \\&\quad+\Qr\lP\mathcal{G}_n,\min_{z\in\textnormal{B}\lP -1+0i, 2tn^{-\nicefrac{11}{10}}\rP} \abs{G_n(z)} \leq tn^{-\nicefrac{1}{2}}\lP \log n\rP^{-\gamma}\rP 
= \textnormal{O}\lP{n^{-1/10}}\rP.
\end{align*}
\subsubsection{\bf Small ball analysis at $\bm e^{\bm{i2\pi x_\alpha}}$} In this part, we are focusing mainly on the complex roots of a Kac polynomial. We remark that the complex roots are more dispersed than the real roots, but they are approaching faster than the real roots to the unit circle. However, the complex roots do not approach extremely fast.

Let $z\in\textnormal{B}( e^{i2\pi x_\alpha},2tn^{-2} \lP\log n\rP^{-1/2-\gamma})$ and assume that $\mathcal{G}_n$ holds.  By Taylor's Theorem we obtain
\begin{eqnarray*}
\abs{G_n(z) - G_n\lP e^{i2\pi x_\alpha}\rP} \leq  \abs{z-e^{i2\pi x_\alpha}}\abs{G'_n\lP e^{i2\pi x_\alpha}\rP} + \abs{R_2(z)},
\end{eqnarray*} where $R_2(z)$ is the error of the Taylor approximation of order 2, and it satisfies
\[
\abs{R_2(z)} \leq  \frac{\lP2tn^{-2}\rP^2}{1-2tn^{-2}} \lC\max_{z\in\C\; :\;\abs{\abs{z}-1} < tn^{-2}} \abs{G_n(z)}\rC  \leq  \frac{4t^2 n^{-5/2}}{1-2tn^{-2}}.
\] 
Then
\begin{align*}
\abs{G_n(z) - G_n\lP e^{i2\pi x_\alpha}\rP} & \leq  2tn^{-2}\lP \log n\rP^{-1/2-\gamma} \| W'_n\| + \frac{4t^2 n^{-5/2}}{1-2tn^{-2}} \\
& \leq  2C_0 tn^{-1/2} \lP\log n\rP^{-\gamma} +\frac{4t^2 n^{-5/2}}{1-2tn^{-2}}.
\end{align*} Hence,
\[
\Qr\lP \mathcal{G}_n, \textnormal{B}_\alpha \rP \leq \Qr\lP\abs{G_n\lP e^{i2\pi x_\alpha}\rP} \leq 2C_4 t n^{-2}\lP \log n\rP^{-\gamma} \rP,
\] 
where $2C_4 = 2C_0 + 4t + 1$.
For proving that $\Qr\lP \mathcal{G}_n, \textnormal{B}_\alpha \rP$ tends to zero as $n\to\infty$, we rewrite the sum $G_n(e^{i2\pi x_\alpha})$ as the product of a matrix by a vector.
This simple 
rewriting allows us to apply lcd techniques for matrices.
To be precise, we define the $2\times n$ matrix $V_\alpha$ as follows
\[
V_\alpha := 
\lC
\begin{array}{cccc}
1 & \cos\lP 2\pi x_\alpha \rP & \ldots & \cos\lP (n-1)2\pi x_\alpha \rP \\
0 & \sin\lP 2\pi x_\alpha \rP & \ldots & \sin\lP (n-1)2\pi x_\alpha \rP
\end{array}
\rC
\] and $X:=\lC\xi_0,\ldots,\xi_{n-1}\rC^T \in \R^n$. Notice that 
\[
V_\alpha X = \lC
\sum_{j=0}^{n-1} \xi_j \cos\lP j2\pi x_\alpha\rP, \sum_{j=0}^{n-1} \xi_j \sin\lP j2\pi x_\alpha\rP \rC^T \in \R^2,
\] which implies
\[
\| V_\alpha X\|_2 = \abs{\sum_{j=0}^{n-1} \xi_j e^{ij2\pi x_\alpha}} = \abs{G_n\lP e^{i2\pi x_\alpha}\rP}.
\]
Let $\Theta = r\lC \cos(\theta), \sin(\theta)\rC^T\in\R^2$, where $r>0$ and $\theta\in\lC0,2\pi\rC$. For fixed $r,\theta$, we have
\[
V_\alpha^T\Theta = r\lC\cos\lP-\theta\rP,\cos\lP2\pi x_\alpha - \theta\rP,\ldots,\cos\lP2\lP n-1\rP \pi x_\alpha - \theta \rP\rC^T.
\]
Note that $\| V_\alpha^T\Theta\|_2\leq r\sqrt{n}$. On the other hand, we have
\[
\det\lP V_\alpha V_\alpha^T \rP = \det\lC
\def\arraystretch{1.5}
\begin{array}{cc}
\sum_{j=0}^{n-1} \cos^2\lP j2\pi x_\alpha \rP & \frac{1}{2}\sum_{j=0}^{n-1} \sin\lP 2\cdot j2\pi x_\alpha \rP \\
\frac{1}{2}\sum_{j=0}^{n-1} \sin\lP 2\cdot j2\pi x_\alpha \rP & \sum_{j=0}^{n-1} \sin^2\lP j2\pi x_\alpha \rP
\end{array}
\rC.
\] 
Now, we are in the setting of 
inequality \eqref{eqn270220191047}.
Recall that $x_\alpha$ satisfies 
\[
n^{-11/10} < \abs{2\pi x_\alpha \hspace{-0.2cm}\mod\pi}<\pi - n^{-11/10}.\] In the following we  distinguish three cases for $x_\alpha$.

\subsubsection{Case 1. $\alpha\in J_1(n,N)$}\label{subsubsec07052020B} Assume that $\gcd\lP\alpha,N\rP \geq n^{1+1/10}\lP\log n\rP^{-\gamma}$. Recall that $N=\lfloor n^2\lP \log n\rP^{1/2+\gamma} \rfloor$. Then we have
\[
\frac{N}{\gcd\lP\alpha,N\rP} \leq \frac{n^2\lP \log n\rP^{1/2+\gamma}}{n^{1+1/10}\lP\log n\rP^{-\gamma}} = n^{1-1/10}\lP \log n\rP^{1/2+2\gamma}.
\] Note that $2\pi x_\alpha$ satisfies $n^{-1}<\abs{2\pi x_\alpha\hspace{-0.2cm}\mod\pi}<\pi - n^{-1}$ for all large $n$. 
By Lemma~3.2 part 1 in \cite{Kon1999}, there exist positive constants $c_5,C_5$ such that
\begin{equation}\label{eqn0403201911759}
c_5 n^2 \leq \det\lP V_\alpha V_\alpha^T\rP \leq C_5 n^2.
\end{equation} Before continue with our arguments, we estimate the number of indexes $\alpha$ where the condition $\gcd\lP \alpha, N\rP \geq n^{1+1/10}\lP\log n\rP^{-\gamma}$ holds. The following lemma provides such estimate.
\begin{lemma}\label{lem040320191725} The number of indices $\alpha$ such that \[\gcd\lP \alpha, N\rP \geq \frac{n^{1+1/10}}{\lP\log n\rP^{\gamma}}\] is at most
\[
n^{1-1/10+\textnormal{o}(1)}\lP\log n\rP^{1/2+2\gamma+\textnormal{o}(1)}.
\]
\end{lemma}
By Proposition~7.4 in \cite{RV2}, the lcd of $V_\alpha$ satisfies $D\lP V_\alpha\rP \geq 1/2$. Thus, by inequalities \eqref{eqn270220191047} and \eqref{eqn0403201911759}, and 
Lemma~\ref{lem040320191725} we obtain
\[
\begin{split}
\def\arraystretch{1.5}
 & \sum_{\alpha\in J_1(n,N)} \Qr\lP \abs{G_n\lP e^{i2\pi x_\alpha} \rP} \leq 2tC_4 n^{-1/2}\lP\log n\rP^{-\gamma} \rP \\ 
& \quad \leq   2n^{1-1/10+\textnormal{o}\lP 1 \rP} \lP \log n\rP^{1/2+2\gamma+\textnormal{o}\lP 1 \rP}\lP \frac{2C^2L^2\lP 2tC_4\rP^2}{\lP c_5 n^2 \rP^{1/2} \lP n^{1/2} \lP\log n\rP^{\gamma}\rP^2} + \frac{2C^2L^2}{\frac{1}{4} \lP c_5 n^2 \rP^{1/2}}\rP \\ 
&\quad = \frac{4C^2L^2 \lP2tC_4\rP^2 \lP\log n\rP^{1/2+\textnormal{o}(1)}}{c_5^{1/2} n^{1+1/10 -\textnormal{o}(1)}} + \frac{4C^2L^2\lP \log n\rP^{1/2+2\gamma+\textnormal{o}(1)}}{\frac{1}{4} c_5^{1/2} n^{1/10-\textnormal{o}(1)}} \\
& \quad \leq C_6 \frac{\lP\log n\rP^{\nicefrac{1}{2}+2\gamma+\textnormal{o}(1)}}{n^{\nicefrac{1}{10}-\textnormal{o}(1)}},
\end{split}
\] where $C_6=4c_5^{-1/2}C^2L^2 \lP\lP2tC_4\rP^2 +4\rP$.
\subsubsection{Case 2. $\alpha\in J_2(n,N)$}\label{subsubsec07052020C} Assume that \[n^{1+1/10} \lP \log n\rP^{-\gamma} \geq \gcd\lP\alpha, N\rP \geq n \lP \log n\rP^{1/2+\gamma}.\] Since $N=\lfloor n^2 \lP\log n\rP^{1/2+\gamma} \rfloor$, we have
\begin{equation}\label{eqn050320191925}
n\geq \frac{N}{\gcd\lP\alpha,N\rP}\geq n^{1-1/10}\lP\log n\rP^{1/2+2\gamma} - \textnormal{o}(1),
\end{equation} where $\textnormal{o}(1)=n^{-1-1/10}\lP \log n \rP^{\gamma}$. 
We observe that $2\pi x_\alpha$ is such that \[n^{-1}\leq \abs{2\pi x_\alpha\hspace{-0.2cm}\mod\pi}\leq \pi - n^{-1}.\] By Lemma~3.2 part 1 in \cite{Kon1999} there exist positive constants $c_5,C_5$ such that 
\[
c_5 n^2 \leq \det\lP V_\alpha V_\alpha^T\rP \leq C_5 n^2.
\]
Also, we observe that $x_\alpha = \frac{\alpha}{N}=\frac{\alpha'}{N'}$ where $\alpha = \alpha' \gcd\lP\alpha,N\rP$ and $N = N'\gcd\lP\alpha, N\rP$. Note that $\gcd\lP\alpha',N'\rP=1$. Since $N'\leq n$, for any $\theta$ we have
\[
\begin{split}
\lL \exp\lP i \lP j2\pi \frac{\alpha'}{N'} - \theta \rP\rP : j = 0,\ldots, N'-1 \rL =  \lL \exp\lP i \lP j2\pi \frac{1}{N'} - \theta \rP\rP : j = 0,\ldots, N'-1 \rL.
\end{split}
\] 
The above observation allows us to assume that $x_\alpha = \nicefrac{1}{N'}$. To apply  inequality \eqref{eqn270220191047}
we need to estimate the lcd.
The following lemma shows an arithmetic property of the values $\cos\lP j2\pi x_\alpha - \theta\rP$ for $j=0,\ldots,N'$ which becomes crucial for estimating the lcd.
\begin{lemma}\label{lem050320191745} 
Fixed $\theta\in[0,2\pi)$ and positive $m\in\Z$. Let $\mathcal{V}$ be a vector in $\R^m$ which entries are $\mathcal{V}_j= r\cos\lP j 2\pi x-\theta \rP$ for $j=0,\ldots,m-1$ with positive integer $r\geq 2$ and $x=\nicefrac{1}{m}$. Then
\[
\mathrm{dist}\lP \mathcal{V},\Z^m\rP \geq \frac{1}{48}\cdot\frac{1}{2\pi x}\quad
\mbox{ whenever } \frac{1}{2r\lP2\pi x\rP}\geq 6.
\] 
\end{lemma}
Since it is needed to analyze \[
V_\alpha^T\Theta = r\lC\cos\lP-\theta\rP,\cos\lP2\pi x_\alpha - \theta\rP,\ldots,\cos\lP2\lP n-1\rP \pi x_\alpha - \theta \rP\rC^T
\] in the definition of  the least common denominator, we can assume without loss of generality that $r$ is a positive integer. In fact, by Proposition 7.4 in \cite{RV2}, we can take $r\geq \nicefrac{1}{2}$. For the case $2>r\geq\nicefrac{1}{2}$, we can replicate the ideas in the proof of Lemma~\ref{lem050320191745} to obtain that $\textnormal{dist}\lP V_\alpha^T\Theta, \Z^n\rP\geq Cn^{1-1/10}$ for some positive constant $C$. If $r\geq 2$, we can use $\lfloor r\rfloor$ instead of $r$ in 
Lemma~\ref{lem050320191745}.

If $r\leq \frac{1}{2\cdot6\cdot2\pi x_\alpha}$, by 
Lemma~\ref{lem050320191745} and expression \eqref{eqn050320191925}, we would obtain
\[
\begin{split}
  \frac{1}{48}\cdot \frac{1}{2\pi} n^{1-1/10} \lP \log n\rP^{1/2+2\gamma} - \textnormal{o}\lP 1\rP & \leq   \frac{1}{48} \cdot \frac{1}{2\pi x_\alpha} 
\leq  \textnormal{dist}\lP V_\alpha^T\Theta, \Z^n\rP \\
& \hspace{-2.5cm}\leq  L\sqrt{\log_+\frac{\| V_\alpha^T \Theta\|_2}{L}} 
 \leq  L\sqrt{\log_+\frac{rn^{1/2}}{L}} \leq L\sqrt{\log_+\frac{n^{3/2} }{L}},
\end{split}
\] which is a contradiction since $L\geq \sqrt{\nicefrac{2}{q}}$ is fixed. Thus, we should have $r > \frac{1}{2\cdot6\cdot2\pi x_\alpha}$ which implies that lcd of $V_\alpha$ satisfies
\[
D\lP V_\alpha \rP > \frac{1}{12}\cdot \frac{1}{2\pi} n^{1-1/10}\lP\log n\rP^{1/2+2\gamma}-\textnormal{o}(1).
\] 
By inequality \eqref{eqn270220191047} we obtain
\[
\begin{split}
 & \sum_{\alpha\in J_2(n,N)} \Qr\lP \abs{G_n\lP e^{i2\pi x_\alpha} \rP} \leq 2tC_4 n^{-1/2}\lP\log n\rP^{-\gamma}\rP \\ 
 & \quad \leq n^2 \lP \log n\rP^{1/2+\gamma} \lP \frac{2C^2L^2 \lP2tC_4\rP^2}{\lP c_5 n^2\rP^{1/2} \lP n^{1/2} \lP\log n \rP^{\gamma}\rP^2} \rP \\ 
 & \quad\quad \;\;+ n^2 \lP \log n\rP^{1/2+\gamma} \lP \frac{2C^2L^2}{\lP c_5 n^2\rP^{1/2} \lP\frac{1}{12}\cdot \frac{1}{2\pi}\cdot n^{1-1/10} \lP \log n\rP^{1/2+2\gamma} - \textnormal{o}\lP 1\rP\rP^2} \rP  \\ 
 & \quad \leq \frac{2C^2L^2 \lP2tC_2\rP^2}{\lP \log n \rP^{\gamma-1/2}} + \frac{2C^2 L^2}{c_5^{1/2}\lP \frac{1}{12}\cdot \frac{1}{2\pi}\rP^2 n^{1-1/5}\lP\log n\rP^{1/2+3\gamma} \lP 1- \textnormal{o}\lP1\rP\rP^2}.\\
 & \quad \leq \frac{C_7}{\lP \log n\rP^{\gamma-\nicefrac{1}{2}}},
\end{split}
\] where $C_7 = 2C^2L^2\lP \lP2tC_2\rP^2 + c_5^{-1/2}\rP$.

\subsubsection{Case 3. $\alpha\in J_3(n,N)$} Assume that $n\lP \log n\rP^{1/2+\gamma} \geq \gcd\lP\alpha, N\rP\geq n^{9/10}\lP\log n\rP^{1/2+\gamma}$. Since that $N=\lfloor n^2 \lP\log n\rP^{1/2+\gamma}\rfloor$, then
\[
n^{11/10} \geq \frac{N}{\gcd\lP\alpha,N\rP} \geq n - \textnormal{o}(1),
\] where $\textnormal{o}(1) = \frac{1}{n\lP\log n\rP^{1/2+\gamma}}$. Note that $2\pi x_\alpha$ satisfies
\[
n^{-11/10} \leq \abs{2\pi x_\alpha\hspace{-0.2cm}\mod\pi} \leq \lP n -\textnormal{o}(1)\rP^{-1}
\] or 
\[
\pi - \lP n -\textnormal{o}(1)\rP^{-1} \leq \abs{2\pi x_\alpha\hspace{-0.2cm}\mod\pi} \leq \pi - n^{-11/10}.
\] By Lemma~3.2 part 2 in \cite{Kon1999}, there exist positive constants $c_5,C_5$ such that
\[
c_5 n^{2-1/5} \leq \det\lP V_\alpha V_\alpha^T\rP \leq C_5n^2.
\]
On the other hand, the number of indexes $\alpha$ which satisfy the condition over $\gcd\lP\alpha,N\rP$ is at most
\[
4N\lP\frac{1}{n-\textnormal{o}(1)} - \frac{1}{n^{1+1/10}}\rP \leq 4n\lP\log n\rP^{1/2+\gamma}\lP \frac{1}{1-\textnormal{o}(1)}-\frac{1}{n^{1/10}}\rP.
\]
In order to use the inequality \eqref{eqn270220191047}, we need to analyze the least common denominator of $V_\alpha$ for this case. In particular, we need to obtain a suitable lower bound for the distance between $V_\alpha^T\Theta$ and $\Z^n$. We use similar ideas using in the proof of Lemma~\ref{lem050320191745}.

As $x_\alpha=\frac{\alpha}{N} = \frac{\alpha'}{N'}$ with $\gcd\lP \alpha',N'\rP = 1$ and $N'\geq n - 1$, then all the points in \[\lL \exp\lP i \lP j 2\pi x_\alpha - \theta \rP\rP : j=0,\ldots, n-1\rL\quad  \textrm{are different}.
\]
Let $r$ be a positive integer and we consider the set of intervals of the form $\lC \frac{m}{r}, \frac{m+1}{r}\rC$ for all $m\in\lC-r,r\rC\cap\Z$. Let $I_m$ and $J_{m}$ be the corresponding arcs on the unit circle whose projection on the horizontal axis is the interval $\lC \frac{m}{r}, \frac{m+1}{r}\rC$. If $r < n$, by the pigeon-hole principle we have that there exists at least one $I_{M}$ (or $J_M$) for some $M\in\lC-r,r\rC\cap\Z$, which contain at least $\nicefrac{n}{2r}$ points $\exp\lP i \lP j 2\pi x_\alpha - \theta \rP\rP$ in it. For each $\cos\lP j 2\pi x_\alpha - \theta\rP \in \lC \frac{M}{r}, \frac{M+1}{r}\rC$, it is defined \[d_j = \min\lL \abs{\cos\lP j 2\pi x_\alpha - \theta\rP - \frac{M}{r}}, \abs{\cos\lP j 2\pi x_\alpha - \theta\rP - \frac{M+1}{r}}\rL.\] Note among the values $d_j$ at most two can be equal and \[\min_{0\;\leq\; l,k\; \leq\; n-1} \lL \abs{l 2\pi x_\alpha - k 2\pi x_\alpha} \rL \geq 2\pi \frac{1}{N'}.\] Observe that for each $0\leq\lambda\leq L$, with $L =\min\lL \lfloor \frac{n}{4\cdot 2r} - \frac{3}{2} \rfloor, \lfloor \frac{N' }{2\cdot 2r\cdot 2\pi} - \frac{1}{2}\rfloor \rL$, there exists at least one $d_j$ such that $d_j\geq \lP 2\lambda + 1\rP 2\pi \frac{1}{N'}$. So, the sum of all $d_j$ is at least
\[
\begin{split}
\sum_{\lambda=0}^{L} \lP 2\lambda + 1\rP 2\pi \frac{1}{N'} & \geq  2\pi\frac{L^2}{N'},
\end{split}
\] and taking $r \leq \lfloor n^{1/4} \rfloor$ it follows that \[2\pi\frac{L^2}{N'} \geq 2\pi\cdot \frac{1}{n^{1+1/10}}\lP\frac{n^{3/4} -\textnormal{o}\lP1\rP}{16\pi}\rP^2 \geq \frac{1}{128\pi}\lP n^{1/4-1/20} - \textnormal{o}\lP1\rP\rP^2.\]
Now, let $v$ be a vector in $\R^n$ whose entries are $v_j=\cos\lP j 2\pi x_\alpha - \theta\rP$ for each $j=0,\ldots,n-1$. If a positive integer $r\leq \lfloor n^{1/4} \rfloor$, by the previous discussion it follows that the vector $rv=(rv_j)_{1\leq j\leq n}$ satisfies
\[
\mathrm{dist}(rv,\Z^n)\geq \frac{1}{128\pi}\lP n^{1/4-1/20} - \textnormal{o}\lP1\rP\rP^2.
\]
Thus, if $r\leq \lfloor n^{1/4} \rfloor$ and taking a fixed $L\geq \sqrt{2/q}$, by the definition of lcd we would deduce that
\[
\begin{split}
\frac{1}{128\pi}\lP n^{1/4-1/20} - \textnormal{o}\lP1\rP\rP^2  \leq   \textnormal{dist}\lP V_\alpha^T \Theta, \Z^n\rP 
& \leq  L\sqrt{ \log_+ \frac{\|V_\alpha^T \Theta\|_2}{L} } \\
& \leq  L\sqrt{\log_+\frac{rn^{1/2}}{L}} \leq L\sqrt{\log_+\frac{n^{3/4}}{L}},
\end{split}
\]
which implies that the lcd of $V_\alpha$ should satisfy $D\lP V_\alpha\rP \geq n^{1/4}$. By \eqref{eqn270220191047}, we obtain
\[
\begin{split}
& \sum_{\alpha\in J_3(n,N)} \Qr\lP \abs{G_n\lP e^{i2\pi x_\alpha} \rP} \leq 2tC_4 n^{-1/2}\lP\log n\rP^{-\gamma}\rP \\
& \quad \leq 4 n\lP\log n\rP^{1/2+\gamma} \lP\frac{1}{1-\textnormal{o}\lP1\rP} - \frac{1}{n^{1/10}}\rP \lP \frac{2C^2L^2\lP 2tC_4\rP^2}{\lP c_5n^{2-1/5}\rP^{1/2}\lP n^{1/2}\lP\log n\rP^{\gamma}\rP^2}\rP \\ 
& \quad\quad +\;4 n\lP\log n\rP^{1/2+\gamma} \lP\frac{1}{1-\textnormal{o}\lP1\rP} - \frac{1}{n^{1/10}}\rP \lP\frac{2C^2L^2}{\lP c_5 n^{2-1/5}\rP^{1/2} \lP n^{1/4} \rP^2} \rP  \\ 
&  \quad=\; 4\lP\frac{1}{1-\textnormal{o}\lP1\rP} - \frac{1}{n^{1/10}}\rP \lP \frac{2C^2 L^2\lP 2t C_4\rP^2}{c_5^{1/2} n^{1-1/10} \lP\log n\rP^{\gamma-1/2}} +  \frac{2 C^2 L^2 \lP\log n\rP^{1/2+\gamma}}{c_5^{1/2} n^{1/2-1/10}}\rP\\
&  \quad\leq C_8\lP\frac{1}{1-\textnormal{o}\lP1\rP}\rP  \frac{\lP \log n\rP^{\nicefrac{1}{2}+\gamma}}{n^{\nicefrac{1}{2}-\nicefrac{1}{10}}},
\end{split}
\] where $C_8=8c_5^{-1/2}C^2 L^2\lP\lP 2t C_4\rP^2 + 1 \rP$.\\

Combining Case 1, Case 2 and Case 3 we obtain
\begin{equation}\label{eq:12}
\sum_{\alpha\in\Lambda} \Qr\lP \mathcal{G}_n, \textnormal{B}_\alpha \rP = \textnormal{O}\left((\log n)^{-\gamma+1/2}\right),\quad
\textrm{where }
\gamma>\nicefrac{1}{2}. 
\end{equation} 
Hence, inequality \eqref{mainestimate} with the help of \eqref{eqn071120190936}, \eqref{eqn091120190937} and \eqref{eq:12}
yields
\[
\Qr\left( \mathcal{M}_n \right) = \textnormal{O}\left((\log n)^{-\gamma+1/2}\right),\quad \textrm{where }
\gamma>\nicefrac{1}{2}.\]
The preceding estimate, inequality \eqref{eq:bound} and relation \eqref{eqn051120191037}
imply Theorem~\ref{thm0108201602}.

\section{Proof of Theorem~\ref{thm01042018}. On the lower bound for the smallest singular value for random circulant matrices}\label{sec060320191036}
\noindent
Let $\rho\in (0,\nicefrac{1}{4})$ be fixed.
We define $x_k=\nicefrac{k}{n}$, $k=0,\ldots,n-1$. 
Note that
\begin{align*}
\Qr\left(s_n(\mathcal{C}_n)\leq  n^{-\rho}\right) & \leq   \sum_{k=0}^{n-1} \Qr \lP \abs{G_n\lP e^{i 2\pi x_k}\rP} \leq n^{-\rho}\rP \\
& \leq  \Qr \lP \abs{G_n\lP 1 \rP} \leq n^{-\rho}\rP + \Qr \lP \abs{G_n\lP -1 \rP} \leq n^{-\rho}\rP \\
& \quad + \sum_{\substack{k=0 \\ k\;:\;  \gcd\lP k, n\rP \;>\; n^{1/2}}}^{n-1}  \Qr \lP \abs{G_n\lP e^{i 2\pi x_k}\rP} \leq n^{-\rho}\rP \\
& \quad +  \sum_{\substack{k=0 \\ k\;:\;  \gcd\lP k, n\rP \;\leq\; n^{1/2}}}^{n-1}  \Qr \lP \abs{G_n\lP e^{i 2\pi x_k}\rP} \leq n^{-\rho}\rP.
\end{align*}
In the sequel, we prove that the right-hand side of the preceding  inequality is  
$\textnormal{O}\lP n^{-2\rho}\rP$. We consider the following three cases.

\paragraph{\textbf{Case 1.}}
The same reasoning using in Section~\ref{subsubsec07052020A} yields
\[
\Qr \lP \abs{G_n\lP 1 \rP} \leq n^{-\rho}\rP + \Qr \lP \abs{G_n\lP -1 \rP} \leq n^{-\rho}\rP = \mbox{O}\lP n^{{-1/2}}\rP.
\]

\paragraph{\textbf{Case 2.} $\gcd\lP k,n\rP > n^{1/2}$.}
By similar reasoning using in the 
first case of the proof of Theorem~\ref{thm0108201602}, 
Section~\ref{subsubsec07052020B}, we deduce 
\begin{align*}
\sum_{\substack{k=0 \\ k\;:\;  \gcd\lP k, n\rP \;>\; n^{1/2}}}^{n-1}  \Qr  \lP \abs{G_n\lP e^{i 2\pi x_k}\rP} \leq n^{-\rho}\rP &\leq  n^{1/2+\textnormal{o}(1)} \lP \frac{2C^2L^2}{c_{5}^{1/2} n^{1+2\rho}} + \frac{2C^2L^2}{ \frac{1}{2}c_{5}^{1/2}n}\rP \\
& \leq   \frac{2C^2L^2}{c_{5}^{1/2}n^{1/2+2\rho-\textnormal{o}(1)}} + \frac{4C^2L^2}{c_5^{1/2} n^{1/2-\textnormal{o}(1)}}\leq  \frac{C_9}{n^{1/2-\textnormal{o}(1)}},
\end{align*} where $C_9 = 4c_5^{-1/2}C^2L^2$. 

\paragraph{\textbf{Case 3.} $\gcd\lP k,n\rP \leq  n^{1/2}$.}
By similar reasoning using in the 
second case of the proof of Theorem~\ref{thm0108201602}, 
Section~\ref{subsubsec07052020C}, we obtain
\begin{align*}
\sum_{\substack{k=0 \\ k\;:\;  \gcd\lP k, n\rP \;\leq\; n^{1/2}}}^{n-1}  \Qr  \lP \abs{G_n\lP e^{i 2\pi x_k}\rP} \leq n^{-\rho}\rP &\leq  n \lP \frac{2C^2L^2}{c_{5}^{1/2} n^{1+2\rho}} + \frac{2C^2L^2}{ c_{5}^{1/2} n \lP \frac{1}{2\cdot 6 \cdot 2\pi}n^{1/2}\rP^2}\rP\\
& \leq  \frac{2C^2L^2}{c_{5}^{1/2}n^{2\rho}} + \frac{1152\pi^2 C^2L^2}{c_5^{1/2} n}\leq  \frac{C_{10}}{n^{2\rho}},
\end{align*} where $C_{10}=c_{5}^{-1/2}C^2L^2\lP 2+1152\pi^2\rP$.\\

\noindent The combination of all the preceding cases yields
$\Qr\left(s_n(\mathcal{C}_n)\leq  n^{-\rho}\right) = \textnormal{O}\lP{n^{-2\rho}}\rP$ for any  $\rho\in (0,\nicefrac{1}{4})$.

\appendix
{
\section{Proofs of 
Lemma~\ref{lem040320191725} and
Lemma~\ref{lem050320191745}}

\begin{proof}[Proof of Lemma~\ref{lem040320191725}]
Write $m:=n^{1+1/10}\lP\log n\rP^{-\gamma}$. Let $T$ be the Euler totient function.  Then we have 
\[
\sum_{\substack{
        \alpha\;:\;
        \gcd(\alpha,N)\; \geq\; m \\
        0\;\leq\; \alpha\;\leq\; N
    } 
}   1 
\leq 
\sum_{\substack{
        d=\lfloor m\rfloor\\
        d\left| N\right.
    } 
}^N   T\lP\frac{N}{d}\rP. 
\]
Notice that $T\lP s\rP\leq s - \sqrt{s}$ for all $s\in\N$. Moreover, if $d(s)$ is the number of divisors of $s$, it is well-known (see Theorem~13.12 in \cite{APO}) that there exists an absolute constant $C>0$ such that
\[
d(s) \leq s^{C\lP\log\log\lP s\rP\rP^{-1}}.
\]
Hence,
\begin{align*} 
\sum_{\substack{
        \alpha \;:\;
        \gcd(\alpha,N)\; \geq\; m \\
        0\;\leq\; \alpha\;\leq\; N
    } 
}   1 
& \leq 
\lP\frac{N}{\lfloor m\rfloor} - \sqrt{\frac{N}{\lfloor m\rfloor}}\rP N^{C\lP\log\log\lP N\rP\rP^{-1}}  \leq  \frac{1}{\lfloor m\rfloor} N^{1+C\lP \log \log N\rP^{-1}}\\
& \leq 2n^{1-1/10+\textnormal{o}(1)}\lP\log n\rP^{1/2+2\gamma+\textnormal{o}(1)},
\end{align*} where $\textnormal{o}(1) = C\lP\log\log\lP N\rP\rP^{-1}$. 
\end{proof}
\begin{proof}[Proof of Lemma~\ref{lem050320191745}]
We define the following sequence \[P=\lL\exp\lP i\lP j2\pi x - \theta\rP\rP: j=0,\ldots,m-1\rL,\] where $i$ is the imaginary unit. Note that $P$ is a set of points on the unit circle which can be seen as vertices of a regular polygon with $m$ sides inscribed in the unit circle. 
Since the arguments of points $\exp\lP i\lP j2\pi x - \theta\rP\rP$ are separated exactly by a distance $2\pi x$, the number of points $\exp\lP i\lP j2\pi x - \theta\rP\rP$ which are in any arc on the unit circle is at least $\frac{l}{2\pi x}-2$, where $l$ is the length of the arc. 

Let $\lC y, y + 3(2\pi x)\rC$ be a subinterval of $[-1,1]$ and consider the arc $A$ on the unit circle whose projection on the horizontal axis is $\lC y, y + 3(2\pi x)\rC$. If the length of the arc $A$ is $l$, then the number of values $\cos\lP j2\pi x - \theta\rP$ which are still in $\lP y, y + 3(2\pi x)\rP$ is at least $\frac{l}{2\pi x}-2\geq \frac{3\lP 2\pi x\rP}{2\pi x} - 2 = 1$ since $l\geq 3\lP 2\pi x\rP$.  

Let $s\in\lC-(r-1),(r-1)\rC\cap\Z$. Note that there exists at least one value \[\cos\lP j2\pi x - \theta\rP\in\lP \frac{s}{r} + 3\lP k-1\rP\lP 2\pi x\rP ,\frac{s}{r} + 3k\lP 2\pi x\rP\rP \subset \lC\frac{s}{r},\frac{s+1}{r}\rC\] for all positive integers $k\leq \frac{1}{3 r \lP 2\pi x\rP}$.
In the sequel, we consider all the values $\cos\lP j2\pi x - \theta\rP\in\lC\frac{s}{r},\frac{s+1}{r}\rC$ and define \[d_j:=\min\lL \abs{\cos\lP j2\pi x - \theta\rP - \frac{s}{r}}, \abs{\cos\lP j2\pi x - \theta\rP - \frac{s+1}{r}}\rL.\] Let $L$ be the biggest integer which satisfies $\lP3\cdot 2 \pi x\rP L \leq \frac{1}{2r}$, or equivalently, $L=\lfloor \frac{1}{2r\lP3\cdot 2 \pi x\rP} \rfloor$. Therefore, the sum of $d_j$ for all $\cos\lP j2\pi x - \theta\rP \in \lC\frac{s}{r},\frac{s+1}{r}\rC$ is at least
\begin{align*}
	\sum_{\lambda=1}^L 2\lambda \lP3\cdot 2\pi x\rP  & =  6\lP 2\pi x\rP \sum_{\lambda=1}^L \lambda \geq  6 \lP 2\pi x \rP \frac{L^2}{2} \\
	& \geq  3\lP 2\pi x \rP \lP \frac{1}{2} \cdot \frac{1}{\lP2r\rP \lP3\cdot 2\pi x\rP} \rP^2 = \frac{1}{12} \cdot \frac{1}{\lP2r\rP^2 \lP2\pi x\rP},
\end{align*} where 
we used
the following inequality 
\[
L \geq \frac{1}{2r\lP3\cdot 2\pi x\rP} - 1 \geq \frac{1}{2} \cdot \frac{1}{2r\lP3\cdot 2\pi x\rP}, 
\] which holds if $\frac{1}{2r\lP2\pi x\rP} \geq 6$.
Let $\sigma_s$ be the sum of $d_j$ for each interval $\lC \frac{s}{r},\frac{s+1}{r}\rC$, $s=-(r-1),\ldots,(r-1)$. As $r\geq 2$, then
\begin{equation*}\label{eqn280221091800}
\sum_{s=-(r-1)}^{r-1} \sigma_s \geq \lP 2r - 2\rP \lP \frac{1}{12}\cdot \frac{1}{\lP2r \rP^2 \lP 2\pi x\rP}\rP \geq \frac{1}{24}\cdot\frac{1}{\lP 2r\rP\lP 2\pi x\rP}.
\end{equation*} By the previous analysis, we have that the distance between the vector $\mathcal{V}\in\R^m$ whose entries are $\mathcal{V}_j = r\cos\lP j2\pi x - \theta\rP$ for $j=0,\ldots,m-1$ with $x=\nicefrac{1}{m}$ to $\Z^m$ is at least
\begin{equation*}\label{eqn280120191810}
r \lP\frac{1}{12}\cdot\frac{1}{\lP 2r \rP\lP2\pi x\rP}\rP = \frac{1}{48}\cdot\frac{1}{2\pi x}, 
\end{equation*} verifying that expression $\frac{1}{2r\lP2\pi x\rP} \geq 6$ is fulfilled.
\end{proof}
}

\section*{Acknowledgments}
The authors would like to thank the constructive and useful suggestions provided by professor Jes\'us L\'opez Estrada.
G. Barrera acknowledges support from a post-doctorate grant held at
Center for Research in Mathematics, (CIMAT, 2015--2016). He would like to express his gratitude to Pacific Ins\-titute for the Mathematical Sciences (PIMS, 2017--2019) for the grant held
at the Department of Mathematical and Statistical Sciences at University of Alberta. He also would like to thank to CIMAT, University of Alberta and  University of Helsinki for all the facilities used along the realization of this manuscript.
P.~Manrique acknowledges support 
from C\'atedras CONACyT-M\'exico for the 
research position held at
Mathematics Institute, Cuernavaca (UNAM, 2017--2020).
He also would like to thank to UNAM for all the facilities used along the realization of this manuscript.

\markboth{}{References}
\bibliographystyle{amsplain}

\begin{thebibliography}{99}

\bibitem{APO} Apostol, T.: Introduction to Analytic Number Theory. {\it Undergraduate Texts in Mathematics}, Springer-Verlag, 1976.

\bibitem{BhaSam2014} Bharucha-Reid, A. \& Sambandham, M.: {\it Random Polynomials: Probability and Mathematical Statistics: a Series of Monographs and Textbooks}. Academic Press, 2014.

\bibitem{BP} Bloch, A. \& P\"olya, G.: On the Roots of Certain Algebraic Equations. \textit{Proceedings of the London Mathematical Society} {\textbf{2}}, no. 1, 1932, 102--114.

\bibitem{BorCha2012} Bordenave, C. \& Chafa\"{i}, D.: Around the Circular Law. {\it Probability Surveys} {\textbf{9}}, 2012, 1--89.

\bibitem{AruRajKou2009} Bose, A., Subhra, R. \& Saha, K.:  Spectral Norm of Circulant-Type Matrices.  {\it Journal of Theoretical Probability} \textbf{24}, no. 2, 2011, 479--516.

\bibitem{chareka} Chareka, P., Chareka, O. \& Kennendy, S.: Locally Sub-Gaussian Random Variables and the Strong Law of the Large Numbers. {\it Atlantic Electronic Journal of Mathematics} \textbf{1}, no. 1, 2006, 75--81.

\bibitem{David2012} Davis, P.: {\it Circulant Matrices}. American Mathematical Society, 2012.

\bibitem{DM} Defant, A. \& Mierczyslaw, L.: Norm Estimates for Random Polynomials on the Scale of Orlicz Spaces. {\it Banach Journal of Mathematical Analysis} \textbf{11}, no. 2, 2017, 335--347.

\bibitem{DNV} Do, Y., Nguyen, O. \& Vu, V.: Roots of Random Polynomials with Coefficients Having Polynomial Growth. {\it Annals of  Probability} \textbf{46}, no. 5, 2018, 2407--2494.

\bibitem{erdos} Erd\"os, P.: Problems and Results on Polynomials and Interpolation. {\it Aspects on Contemporary Complex Analysis} \textbf{1980}, 383--391.

\bibitem{HP} Huhtanen, M. \& Per\"am\"aki,  A.: Factoring Matrices into the Product of Circulant and Diagonal Matrices. {\it Journal of Fourier Analysis and Applications} \textbf{21}, no. 5, 2015, 1018--1033.

\bibitem{IbrZap2013} Ibragimov, I. \& Zaporozhets, D.: On Distribution of Zeros of Random Polynomials in Complex Plane, {\it Prokhorov and Contemporary Probability Theory, Proceedings in Mathematics and Statistics} 33, Editors: A. Shiryaev, S. Varadhan, E. Presman, Springer, 2013, 303--323.

\bibitem{kac1943} Kac, M.: A Correction to On the Average Number of Real Roots of a Random Algebraic Equation. {\it Bulletin of the American Mathematical Society} \textbf{49}, no. 1, 1943, 314--320.

\bibitem{Kah1985} Kahane, J.:  Some Random Series of Functions. {\it Cambridge University Press}. Cambridge, Second Edition, 1985.

\bibitem{Kar1997} Karapetyan, A.: On Minimum Modulus of Trigonometric Polynomials with Random Coefficients. {\it Mathematical Notes} \textbf{61}, 1997, 369--373.

\bibitem{Kar1998} Karapetyan, A.: The Values of Stochastic Polynomials in a Neighbourhood of the Unit Circle. {\it Mathematical Notes} 
 \textbf{63}, 1998, 127--130.

\bibitem{KHA} Khattree, R.: Multidimensional Statistical Analysis and Theory of Random Matrices. {\it Proceedings of the Sixth Eugene Lukacs Symposium}, Bowling Green, Ohio, USA, 1996, 29--30.

\bibitem{Kon1994} Konyagin, S.: Minimum of the Absolute Value of Random Trigonometric Polynomials with Coefficients $\pm 1$. {\it Mathematical Notes} \textbf{56}, 1994, 931--947.

\bibitem{Kon1999} Konyagin, S. \& Schlag, W.: Lower Bounds for the Absolute Value of Random Polynomials on a Neighbourhood of the Unit Circle. {\it Transactions of the American Mathematical Society} \textbf{351}, 1999, 4963--4980.

\bibitem{larry} Larry, A. \&  Vanderbei, R.: The Complex Zeros of Random Polynomials. {\it Transactions of the American Mathematical Society} \textbf{347}, no. 11,  1995, 4365--4384.

\bibitem{Lub2016} Lubinsky, D., Pritsker, I., \& Xie, X.: Expected Number of Real Zeros for Random Linear Combinations of Orthogonal Polynomials. {\it Proceedings of the American Mathematical Society} \textbf{144}, no. 4, 2016, 1631--1642.

\bibitem{Mar2009} Meckes, M.: Some Results on Random Circulant Matrices. {\it High Dimensional Probability: The Luminy V}, Institute of Mathematical Statistics, 2009, 213--223.

\bibitem{Mezin} Mezincescu, G., Bessis, D., Fournier, J., Mantica, G. \& Aaron, F.:   Distribution of Roots of Random Real Generalized Polynomials. {\it Journal of Statistical Physics} \textbf{86}, no. 3--4,  1997, 675--705.

\bibitem{Nguyen2016} Nguyen,  H., Nguyen, O., \& Vu, V.: On the Number of Real Roots of Random Polynomial. {\it Communications in Contemporary Mathematics} \textbf{18}, 
no. 4, 2016.

\bibitem{RahSch2002} Rahman Q. \& Schmeisser, G.:  Analytic Theory of Polynomials: Critical Points, Zeros and Extremal Properties. {\it Oxford Science Publications} 2002.

\bibitem{RAU} Rauhut, H.: Circulant and Toeplitz Matrices in Compressed Sensing. {\it Proc. SPARS'09}, Saint-Malo, France, 2009.

\bibitem{RV} Rudelson, M. \& Vershynin, R.: Non-Asymptotic Theory of Random Matrices: Extreme Singular Values. {\it Proceedings of the International Congress of Mathematicians Hyderabad Volume III, India, Editor: Rajendra Bhatia}, 2010, 1576--1602.

\bibitem{RV1} Rudelson, M. \& Vershynin, R.: The Littlewood--Offord Problem and Invertibility of Random Matrices.
\textit{Advances in Mathematics} \textbf{218}, 2008, 600--633.

\bibitem{RV2} Rudelson, M. \& Vershynin, R.: No-gaps Delocalization for General Random Matrices. {\it Geometric and Functional Analysis} \textbf{26}, no. 6, 2016, 1716--1776.

\bibitem{sazy} Salem, R. \& Zygmund A.: Some Properties of Trigonometric Series Whose Terms Have Random Signs. {\it Acta Mathematica} \textbf{91}, 1954, 245--301.

\bibitem{SenVir2013} Sen, A. \& Vir\'ag, B.: The Top Eigenvalue of the Random Toeplitz Matrix and the Sine Kernel. {\it Annals of Probability} \textbf{41}, no. 6, 2013, 4050--4079.

\bibitem{tikomirov2016} Tikhomirov, K.: The Smallest Singular Value of Random Rectangular Matrices with no Moment Assumptions on Entries. {\it Israel Journal of Mathematics} {\bf{212}}, 2016, 289--314.

\bibitem{weber} Weber, M.: On a Stronger Form of Salem--Zygmund's Inequality for Random Trigonometric Sums with Examples. {\it Periodica Mathematica Hungarica} {\bf{52}}, no. 2, 2006, 73--104.


\bibitem{roman} Vershynin, R.: High-Dimensional Probability. An Introduction with Applications in Data Science. {\it Cambridge Series in Statistical and Probabilistic
Mathematics} {\bf{47}}. Cambridge University Press, Cambridge, 2018.
\end{thebibliography}

\end{document}